\documentclass[preprint,12pt,3p,authoryear]{elsarticle}

\usepackage{graphicx}
\usepackage{bm}
\usepackage{booktabs}
\usepackage{array}
\usepackage{tabularx}
\usepackage{url}
\usepackage{caption}
\usepackage{subcaption}
\usepackage{amsmath}
\usepackage{amssymb}
\usepackage{amsthm}
\usepackage{float}
\usepackage[section]{placeins}
\usepackage[T1]{fontenc}
\usepackage{newtxtext,newtxmath}
\theoremstyle{plain}
\newtheorem{proposition}{Proposition}
\newtheorem{corollary}[proposition]{Corollary}

\newcommand{\jump}[1]{[\![#1]\!]}
\newcommand{\Ylog}{\bm{Y}^{\log}}
\newcommand{\Hlog}{\bm{H}^{\log}}

\journal{}

\begin{document}

\begin{frontmatter}

\title{Layer-wise discrete certificates in multicomponent reacting-flow
discretizations}

% \author[a,b]{Jingchao Zhang}
\author[a,b]{Jingchao Zhang\corref{cor1}}
\ead{jingchaozhang@buaa.edu.cn}
\cortext[cor1]{Corresponding author}
\affiliation[a]{organization={Hangzhou International Innovation Institute,
Beihang University}, city={Hangzhou}, postcode={311115}, country={China}}
\affiliation[b]{organization={School of Aeronautic Science and Engineering,
Beihang University}, city={Beijing}, postcode={100191}, country={China}}

\begin{abstract}
Multicomponent reacting-flow discretizations must conserve total species mass, preserve positive mass fractions, and dissipate a convex mixing potential. These properties are controlled at different stages of a solver, however, and a certificate established at one stage need not survive the next. We study this compatibility problem for mixture-averaged diffusion, non-orthogonal finite-volume reconstruction, and positivity-preserving Patankar time integration. The analysis is isothermal, uses prescribed density, and uses synthetic detailed-balance reaction networks for the chemistry tests. At the face level, the entropy-variable jump is shown to equal a log-mean Gibbs Hessian applied exactly to the mass-fraction jump. Requiring conservation and entropy neutrality of the mass-correction term on the simplex tangent space then yields a unique componentwise-positive face composition and an exact two-point entropy certificate. The arithmetic-weight face criterion is also satisfied by production neutral-air/carbon data and a quasi-neutral ambipolar-air control, but a face-level certificate does not automatically extend to subsequent discretization layers. A weighted least-squares non-orthogonal correction can, in the tested prototype, convert a dissipative two-point operator into an anti-dissipative multipoint operator. Likewise, Patankar rescaling can preserve positivity and conservation while increasing the free energy; an explicit finite-step rational counterexample establishes this failure independently of ensemble statistics. A fully discrete two-cell example further shows that a positive face entropy margin does not guarantee positivity of the fully discrete update. Thus conservation, positivity, and free-energy dissipation are layer-specific properties and require distinct discrete certificates.
\end{abstract}

\begin{keyword}
multicomponent diffusion \sep discrete entropy \sep positivity \sep
finite volume \sep Patankar schemes \sep structure-preserving discretization
\end{keyword}

\end{frontmatter}

\section{Introduction}
\label{sec:intro}

High-enthalpy multicomponent flows involve strong species gradients, molecular
diffusion, and chemical reactions evolving over different time scales. Numerical
discretizations of such flows must preserve several basic structural properties.
Total species mass must remain conserved, mass fractions must stay positive, and
an appropriate thermodynamic potential should not increase under dissipative
transport and detailed-balance chemistry. Violating these properties can lead to
composition drift, nonphysical species states, or artificial growth of the
discrete free energy. The difficulty is that these constraints are not imposed
at a single stage of the numerical method. Conservation is affected by the
construction of the diffusive flux, positivity depends on both spatial transport
and time integration, and thermodynamic dissipation depends on how the transport
closure, reconstructed gradients, and time update interact. A discretization may
therefore satisfy each requirement in isolation without preserving all of them
after its individual components are combined. The central question is thus not
only how each property can be enforced, but also whether a discrete guarantee
established at one stage of the solver remains valid after the next stage is
introduced.

Considerable progress has been made on these properties separately.
Entropy-conservative and entropy-stable two-point discretizations are available
for multicomponent compressible Euler systems
\citep{tadmor2003,gouasmi2020,renac2021} and have been extended to thermally
perfect gases \citep{aiello2026}, with convergence analyses developed for
structure-preserving multicomponent schemes \citep{agnihotri2026}. For
diffusive systems, entropy production has been established for regularized
gas-mixture models \citep{zlotnik2022,zlotnik2023}, while fully discrete
cross-diffusion schemes can preserve entropy, positivity, and volume-filling
constraints when the face mobility is constructed together with an appropriate
vector-valued discrete chain rule \citep{jungel2023,cances2023}. Positivity
preservation has also been studied extensively for stiff
production--destruction systems. Modified Patankar--Runge--Kutta schemes provide
unconditional positivity together with second- or third-order accuracy
\citep{burchard2003,kopecz2018apnum,kopecz2018bit}, and subsequent work has
clarified their stability, convergence, and relaxation properties
\citep{izgin2022lyapunov,izgin2023mpdec,bender2025,izgin2026relaxation}.
These developments provide strong guarantees within their respective settings.
What remains less clear is whether those guarantees remain compatible when the
spatial and temporal components are assembled into a multicomponent
reacting-flow discretization. This compatibility question is especially
important in reacting-flow solvers, where the diffusion closure, mesh
reconstruction and stiff-chemistry integrator are typically designed
independently.

This question becomes particularly important for three commonly used numerical
ingredients. First, the Hirschfelder--Curtiss mixture-averaged diffusion model
\citep{hirschfelder1954,ern1994,giovangigli1999,cantera2026,comsol2026}
enforces zero total diffusive mass flux through a mass-correction term, but the
resulting diffusion matrix is generally non-symmetric. It therefore does not
fall directly within the symmetric positive-semidefinite framework commonly
used in entropy-stable viscous discretizations \citep{hughes1986}. In addition,
the thermodynamic structure of a two-point face flux depends on how the face
composition is defined, so mass conservation alone does not determine its
entropy behaviour. Second, on non-orthogonal finite-volume meshes, the
two-point normal flux is commonly supplemented by a reconstructed tangential
correction. Weighted least-squares reconstruction improves geometric
consistency, but it also enlarges the stencil from a two-cell interaction to a
multipoint operator, so a face-wise entropy argument cannot be assumed to carry
over unchanged. Third, Patankar-type time integrators preserve positivity by
rescaling production and destruction terms, whereas free-energy decay imposes a
different condition on the discrete update. Stability of a Patankar fixed point
therefore does not by itself imply monotone decay of a thermodynamic potential
along a time-discrete trajectory
\citep{torlo2022,izgin2022lyapunov,izgin2023mpdec}. These three issues have the
same underlying form: a structural property established for one discrete
component need not be preserved when the next component is added.

The present work complements these analyses of order, stability and fixed
points by evaluating dissipation of a specified convex detailed-balance free
energy, constructing an exact rational finite-step counterexample, and
examining six representative third-order members.

The present work examines this problem by separating the discretization into
three layers rather than assigning structural properties to the solver as a
whole. At the face level, we analyse the relation between mass-fraction jumps,
entropy-variable jumps, the mixture-averaged diffusion matrix, and the
mass-correction term on the composition simplex. We then examine how this
face-level structure changes when a cell-centred weighted least-squares
correction is introduced on non-orthogonal meshes. The temporal layer is studied
separately for detailed-balance production--destruction systems advanced with
positivity-preserving Patankar methods, so that conservation, positivity, and
thermodynamic dissipation can be assessed independently. Finally, the spatial
and temporal results are brought together at the fully discrete level to
determine which structural guarantees can be combined and which require an
additional condition. This layer-wise formulation keeps the assumptions and
scope of each argument explicit and provides a common framework for analysing
how conservation, positivity, and dissipation interact across the different
parts of a reacting-flow discretization.

\section{Governing structure and discrete observables}
\label{sec:structure}

\subsection{Composition state and mixing potential}

Consider $N$ species with molar masses $M_s$ and mass fractions $\bm{Y}$ on the
simplex $\sum_s Y_s=1$, $Y_s>0$, and write $b_s=1/M_s$, $Z=\sum_s Y_s/M_s$ and
$X_s=(Y_s/M_s)/Z$. At fixed temperature we use the isothermal ideal-mixture
Gibbs composition potential; the prescribed density is a coefficient rather
than an equation-of-state variable. The $\ln Z$ term is the Gibbs-form
composition contribution and gives $\ker H=\operatorname{span}\{Y\}$. A
fixed-$(T,\rho)$ Helmholtz potential would omit this term and would have a
diagonal composition Hessian; the face construction in Section~\ref{sec:layer1}
is therefore specific to the Gibbs form used here. The potential $\varphi$ has gradient and
Hessian
\begin{equation}
w_s=\frac{\partial\varphi}{\partial Y_s}=a_s+\frac{R_u}{M_s}\ln X_s,
\qquad
\bm{H}=\frac{\partial\bm{w}}{\partial\bm{Y}}
=R_u\!\left[\operatorname{diag}\!\Big(\tfrac{1}{M_sY_s}\Big)
-\frac{\bm{b}\bm{b}^{\mathsf T}}{Z}\right],
\label{eq:potential}
\end{equation}
with $\bm{H}\bm{Y}=0$ and $\bm{H}$ positive definite on the tangent space
$T=\{\bm{z}:\bm{1}^{\mathsf T}\bm{z}=0\}$. The vector $\bm{w}$ plays the role of
the entropy variables. Here $R_u$ is the universal gas constant; numerically
we use $R_u=8.314462618\,{\rm J\,mol^{-1}\,K^{-1}}$, and $a_s$
is the composition-independent standard-state contribution to $w_s$ at
the prescribed thermodynamic state; the source implementation writes it as
$a_s=g_s^0/M_s$, where $g_s^0$ is the standard molar potential divided by
temperature. Throughout, $\operatorname{sym}(B)=(B+B^{\mathsf T})/2$.

Throughout, ``entropy dissipation'' means non-increase of the convex mixing
potential $\varphi$ defined in \eqref{eq:potential}; the terms ``mixing
potential'' and ``free energy'' refer to this same quantity, with the sign
convention used in \eqref{eq:potential_total}.

To isolate the composition-space mechanisms, temperature is fixed and density
is treated as a prescribed positive, time-independent coefficient in the
species-transport problem considered here. Fourier conduction
and the Soret and Dufour cross effects are absent, as is pressure diffusion.
The reaction source is isothermal and the gas--surface channel is not
represented. What remains is the smallest system in which the three constraints
of Section~\ref{sec:intro} can compete: species transport by mixture-averaged
diffusion, a stiff detailed-balance source, and a mesh that need not be
orthogonal.

Equation~\eqref{eq:semidiscrete} uses the reduced mixture-averaged closure
$\bm J_s=-\rho d_s\nabla Y_s+Y_s\sum_k\rho d_k\nabla Y_k$.
It is a deliberate composition-space model: the effective coefficients
$d_s$ are supplied by the Chapman--Enskog transport controls used below, while
barodiffusion, thermal diffusion and the explicit $\nabla\overline M$ term of
more complete mole-fraction formulations are outside the present scope.

\subsection{The semi-discrete system and its three discrete choices}
\label{sec:structure:system}

On a mesh of cells $i$ with volumes $V_i$ and faces $f$ of area $|f|$,
let $s_{if}=\pm1$ denote the orientation of an interior face relative to cell
$i$. The semi-discrete system studied throughout is
\begin{equation}
V_i\,\rho_i\,\frac{\mathrm{d}\bm{Y}_i}{\mathrm{d}t}
=-\sum_{f\in\partial i}s_{if}|f|\,\bm{J}_f
+V_i\,\bm{\omega}(\bm{Y}_i),
\qquad
\bm{J}_f=-\rho_f\,\bm{A}_f\,\bm{g}_f .
\label{eq:semidiscrete}
\end{equation}
Here $\rho_i>0$ is the prescribed, time-independent density in cell $i$ and
$\rho_f>0$ is the prescribed face density. The vector $\bm{J}_f$ is oriented
from the left cell to the right cell; the sign $s_{if}$ supplies the outward
orientation in the cell balance. The vector $\bm{\omega}$ is the source rate
per unit volume. The vector $\bm d_f=(d_{1,f},\ldots,d_{N,f})^{\mathsf T}$
collects the face mixture-averaged diffusivities.
Three quantities in \eqref{eq:semidiscrete} are not supplied by the continuous
equations and must be chosen when the scheme is built. They are the three
layers of this paper.

We call a discrete construction certified for a property when an explicit
identity, inequality or invariant proves that property for the stated operator
and hypotheses. Finite-population results are reported as observed rather than
certified.

\begin{figure}[htbp]
\centering
\includegraphics[width=0.98\textwidth]{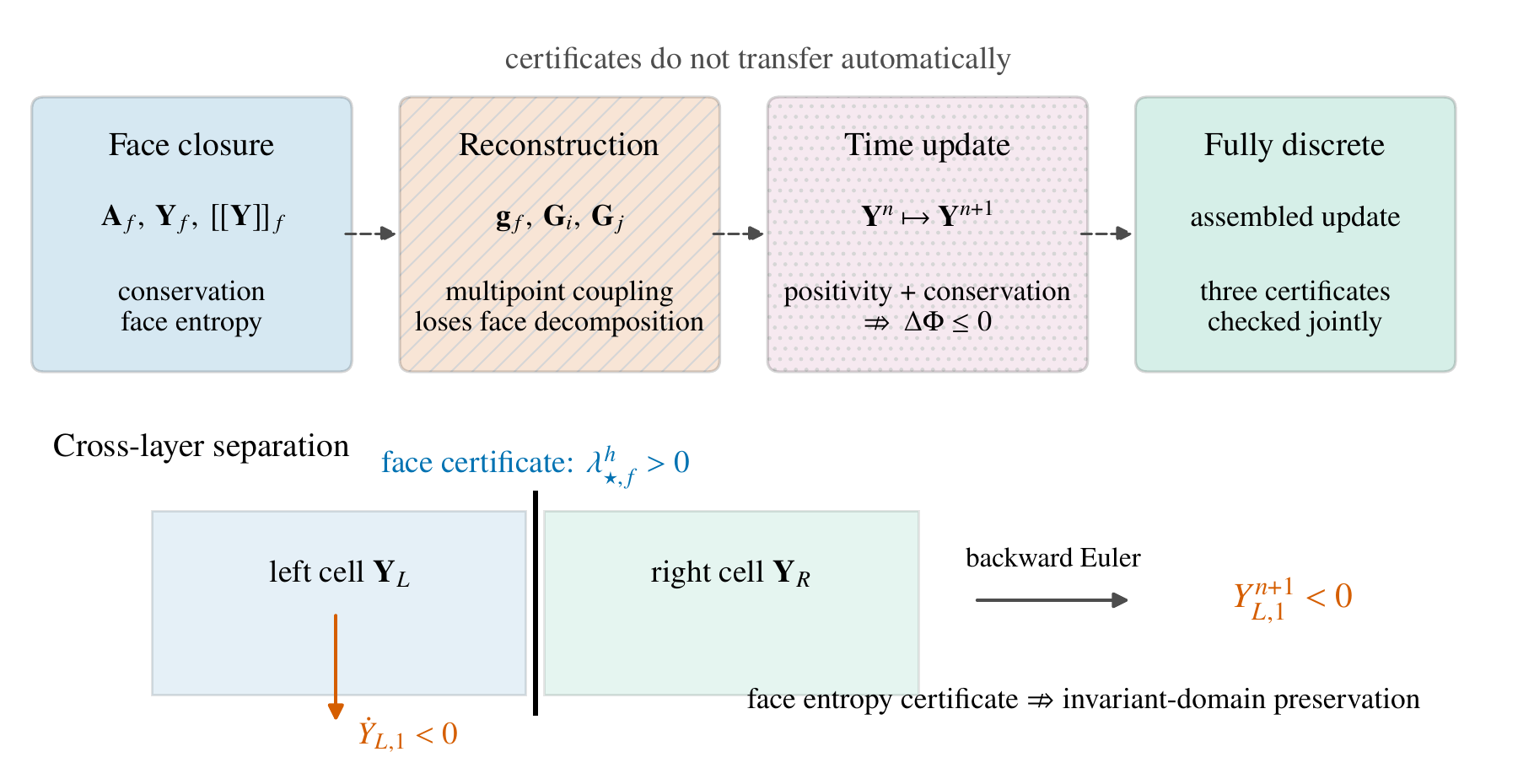}
% [图源] figures/scripts/layerwise_certificate_map.py | 数据: analytical schematic, not simulation data
\caption{Conceptual map of the three discrete layers and their certificates.
The dashed links are labelled ``certificates do not transfer automatically''
and indicate that a certificate established for one algebraic object is not
inherited automatically by the next. The lower schematic is the
two-cell fully discrete separation used in Section~\ref{sec:baseline}: a
positive face entropy margin can coexist with an outward spatial direction and
a negative component after a backward-Euler step. This is an analytical map,
not a numerical-data figure.}
\label{fig:certificate-map}
\end{figure}

Figure~\ref{fig:certificate-map} summarizes the layer-wise logic before the
individual results are developed.

The first is the face matrix together with the face composition,
\begin{equation}
\bm{A}_f=\operatorname{diag}(\bm{d}_f)-\bm{Y}_f\bm{d}_f^{\mathsf T},
\qquad \bm{1}^{\mathsf T}\bm{Y}_f=1 ,
\label{eq:closure}
\end{equation}
in which the rank-one term is the mass-correction flux that forces
$\bm{1}^{\mathsf T}\bm{J}_f=0$. The face weight is a discretization choice not
specified by the continuous closure. Section~\ref{sec:layer1} treats this layer.

For every oriented interior face, the left cell is $i=L$ and the right cell
is $j=R$. Its unit normal $\bm n_f$ points from left to right, and
$\jump{q}_f=q_R-q_L$. With cell centroids $\bm x_i,\bm x_j$, set
$\bm r_f=\bm x_j-\bm x_i$,
$d_{n,f}=\bm r_f\cdot\bm n_f>0$ and
$\bm t_f=\bm r_f-d_{n,f}\bm n_f$.

\begin{figure}[htbp]
\centering
\includegraphics[width=0.86\textwidth]{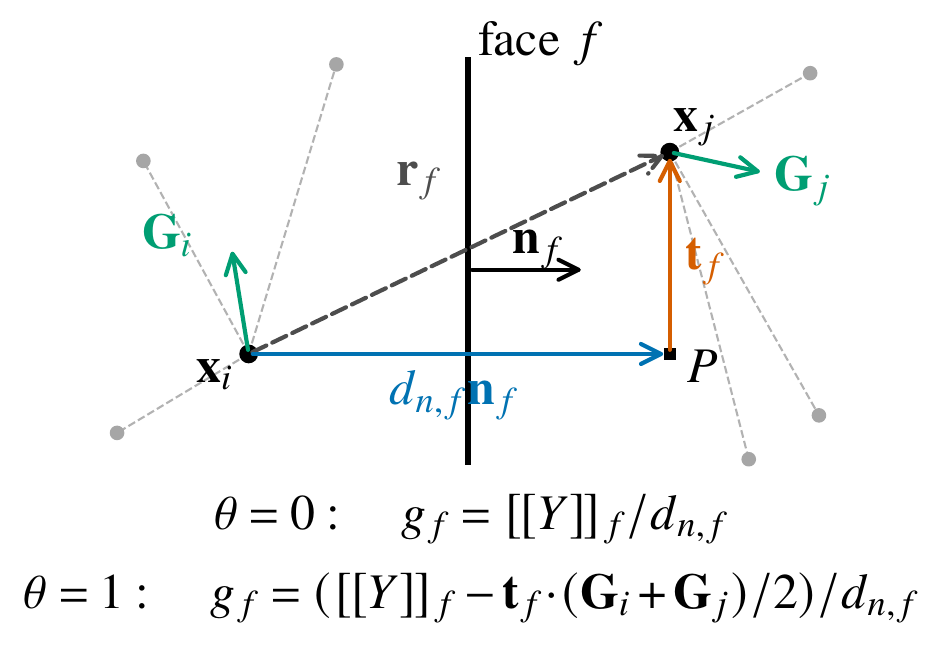}
% [图源] figures/scripts/face_geometry_schematic.py | 数据: analytic illustrative geometry, not simulation data
\caption{Non-orthogonal face geometry, shown schematically rather than as a
computed mesh. The centroid displacement is decomposed as
$\bm r_f=d_{n,f}\bm n_f+\bm t_f$; $P$ is the normal projection of
$\bm x_j$ onto the line through $\bm x_i$. Grey dashed links indicate
illustrative WLS neighbours, not additional transport faces.
At $\theta=0$ the normal derivative uses only the composition jump divided
by $d_{n,f}$. At $\theta=1$ it subtracts the tangential-gradient contribution
in \eqref{eq:reconstruction}. Neither choice changes the physical face location.}
\label{fig:geometry}
\end{figure}

The second is the reconstructed normal gradient $\bm{g}_f$. On a polygonal mesh
the segment joining two cell centroids is not parallel to the face normal.
Writing $d_{n,f}$ for the normal separation and $\bm{t}_f$ for the tangential
offset, and using cell gradients $\bm{G}_i$ from inverse-square weighted least
squares (WLS), where $\bm G_i\in\mathbb R^{N\times d}$ collects the
reconstructed species gradients and row $s$ is $(\nabla Y_s)_i^{\mathsf T}$,
\begin{equation}
\bm{g}_f=\frac{\jump{\bm{Y}}
-\theta\,\tfrac12(\bm{G}_i+\bm{G}_j)\bm{t}_f}{d_{n,f}},
\qquad
\eta_h=\max_f\frac{\|\bm{t}_f\|}{d_{n,f}} .
\label{eq:reconstruction}
\end{equation}
The scalar $\theta\in[0,1]$ scales the tangential correction: $\theta=0$ leaves
the projected two-point operator and $\theta=1$ the full correction
\citep{jasak1996,mavriplis2003}. The quantity $\eta_h$ measures the relative
non-orthogonality of the mesh. Section~\ref{sec:layer2} treats this layer.
Figure~\ref{fig:geometry} shows the centroid decomposition and the
separate roles of the two-point jump and the WLS correction.

The third is the time update that advances \eqref{eq:semidiscrete}, treated in
Section~\ref{sec:layer3}.

Against these choices stand three observables, evaluated on the discrete
solution and not on its continuous limit: the simplex residual
$|\bm{1}^{\mathsf T}\bm{Y}_i-1|$, the smallest mass fraction $\min_{i,s}Y_{s,i}$,
and the free energy
\begin{equation}
\Phi=\sum_i V_i\,\rho_i\,\varphi(\bm{Y}_i).
\label{eq:potential_total}
\end{equation}
A fully discrete step uses $\Delta\Phi:=\Phi^{n+1}-\Phi^n$; in a
single-cell or source-only calculation we write $\Phi_0:=\Phi(\bm Y^0)$.
A scheme is called conservative, positive or dissipative according to whether
the first stays at machine zero, the second stays above zero, and the third
does not increase.

\subsection{Detailed-balance sources and Patankar updates}

The source is written as a conservative production--destruction system (PDS),
\begin{equation}
\rho\,\frac{\mathrm{d}Y_s}{\mathrm{d}t}=\sum_{r}\big(p_{sr}(\bm{Y})-p_{rs}(\bm{Y})\big),
\qquad p_{sr}\ge0 ,
\label{eq:pds}
\end{equation}
built from detailed balance so that the continuous entropy production
$\sigma=-\bm{w}\cdot\bm{\omega}\ge0$ holds term by term. The antisymmetric form
of \eqref{eq:pds} makes $\sum_s\rho\,\mathrm{d}Y_s/\mathrm{d}t=0$ structural.
For each elementary reaction $r$, let
$f_{sr}=M_s\nu_{sr}q_r$ be its signed mass rate, with positive entries
classified as gains and negative entries as losses. The PDS splitting used in
all Patankar calculations assigns
$p_{sq}^{(r)}=\max(f_{sr},0)\max(-f_{qr},0)/\sum_k\max(f_{kr},0)$
and sums this expression over reactions; the destruction matrix is its
transpose. Thus mass released by all losing species is distributed
proportionally over the gaining species, exactly as in the archived source.

A Patankar-type update multiplies each destruction term by the ratio of the
unknown to a known value \citep{patankar1980,burchard2003}, which for the
first-order member gives
\begin{equation}
\rho\,\frac{Y_s^{n+1}-Y_s^{n}}{\Delta t}
=\sum_{r}\left(p_{sr}(\bm{Y}^n)\frac{Y_r^{n+1}}{Y_r^{n}}
-p_{rs}(\bm{Y}^n)\frac{Y_s^{n+1}}{Y_s^{n}}\right).
\label{eq:mpe}
\end{equation}
The construction targets two of the three observables: \eqref{eq:mpe} is
unconditionally positive and exactly conservative, and it costs one linear
solve. Nothing in it constrains $\Phi$, and Section~\ref{sec:layer3} shows that
the missing inequality does not follow from the other two.

\section{Analytical results and numerical evidence}
\label{sec:verification}

We distinguish analytical results from finite numerical evidence throughout.
Exact identities and propositions follow from the stated
hypotheses, with symbolic algebra used as an independent cross-check where
applicable. Numerical counterexamples and finite-population observations are
additionally checked by independent reassembly or external implementations as
noted. Table~\ref{tab:props} summarizes the analytical results and hypotheses.
Appendix B gives the population definitions, numerical tolerances and
resampling procedures.

\begin{table}[htbp]
\centering
\caption{Analytical results and their hypotheses. All rows assume the standing
conditions of Section~\ref{sec:structure}; additional hypotheses are listed
explicitly. All propositions and
corollaries in this table are proved from the stated hypotheses. Symbolic and
numerical checks are independent verification controls and do not form part of
the analytical hypotheses.}
\label{tab:props}
\small
\begin{tabularx}{\textwidth}{@{}l>{\raggedright\arraybackslash}X>{\raggedright\arraybackslash}Xl@{}}
\toprule
 & Content & Hypotheses & Status \\
\midrule
Prop.~\ref{prop:identity} & $\jump{\bm{w}}=\Hlog\jump{\bm{Y}}$
  & $\bm{Y}>0$ & proved \\
Cor.~\ref{cor:telescope} & face-wise telescoping
  & closed mesh, $\theta=0$ & proved \\
Prop.~\ref{prop:kernel} & $\ker\Hlog$ non-trivial iff $\zeta=1$
  & $\bm{Y}>0$ & proved \\
Prop.~\ref{prop:tangentweight} & conservative tangent-neutral weight
  & $\bm{Y}_{L,R}>0$ & proved \\
Cor.~\ref{cor:tangent_reduction} & tangent-space criterion reduces to Fickian part
  & $\bm{Y}_f=\bm{Y}_f^\star$, $\bm{z}\in T$ & proved \\
Prop.~\ref{prop:shift} & generalised margin shifts by $\delta$ exactly
  & $\bm{1}^{\mathsf T}\bm{Y}_f=1$ & proved \\
Prop.~\ref{prop:ascent} & ascent resists scalar rescaling
  & $\Phi$ convex and differentiable, ray in domain & proved \\
\bottomrule
\end{tabularx}
\end{table}

\section{Layer 1: the two-point face certificate}
\label{sec:layer1}

\subsection{An exact identity and an exact criterion}

\begin{proposition}[Log-mean representation]
\label{prop:identity}
Define the log-mean Gibbs Hessian
\begin{equation}
\Hlog(\bm{Y}_L,\bm{Y}_R)
=R_u\!\left[\operatorname{diag}\!\Big(\tfrac{1}{M_sY_s^{\log}}\Big)
-\frac{\bm{b}\bm{b}^{\mathsf T}}{Z^{\log}}\right],
\label{eq:hlog}
\end{equation}
with $Y_s^{\log}=\mathrm{logmean}(Y_{s,L},Y_{s,R})$ and
$Z^{\log}=\mathrm{logmean}(Z_L,Z_R)$. Then
$\jump{\bm{w}}=\Hlog\jump{\bm{Y}}$ holds pointwise and exactly.
\end{proposition}

The proof is one line from $\jump{\ln a}=\jump{a}/a^{\log}$ and is given below.
The identity carries no truncation error, which is
what makes the following criterion exact rather than asymptotic.

\begin{corollary}[Discrete entropy telescoping]
\label{cor:telescope}
For the transport part of \eqref{eq:semidiscrete} on a closed mesh (periodic
or zero-flux boundary faces) with $\theta=0$,
\begin{equation}
\frac{\mathrm{d}\Phi}{\mathrm{d}t}
=-\sum_f\frac{\rho_f|f|}{d_{n,f}}\,\jump{\bm{Y}}^{\mathsf T}\bm{K}^h_f\,\jump{\bm{Y}},
\qquad
\bm{K}^h_f=\operatorname{sym}\!\left[\Hlog\bm{A}_f\right],
\label{eq:telescope}
\end{equation}
for arbitrary positive face areas and normal separations. Thus a uniform
semi-discrete face certificate is
$\bm{Q}^{\mathsf T}\bm{K}^h_f\bm{Q}\succeq0$, where the columns of
$\bm{Q}\in\mathbb{R}^{N\times(N-1)}$ are an orthonormal basis of $T$.
\end{corollary}

Let $\lambda_{\star,f}^h$ be the smallest generalised eigenvalue of
$(\bm{Q}^{\mathsf T}\bm{K}^h\bm{Q},\ \bm{Q}^{\mathsf T}\Hlog\bm{Q})$.

Equation~\eqref{eq:telescope} is an algebraic identity, not an approximation.
The Hessian is evaluated at the log-mean composition, so a criterion evaluated
at either nodal state is not the discrete face criterion. Beyond the floating-point
regression, Propositions~\ref{prop:identity}, \ref{prop:kernel},
\ref{prop:tangentweight} and \ref{prop:shift} were re-derived in exact symbolic
algebra with free symbols for $N=2,3,4$. The new face weight was additionally
checked on $35\,000$ positive faces with $N=2,\ldots,8$.

\subsection{A conservative tangent-neutral face weight}

Splitting $\bm{K}^h_f$ separates the Fickian part from the correction part,
\begin{equation}
\bm{K}^h_f
=\operatorname{sym}\!\left[\Hlog\operatorname{diag}(\bm{d}_f)\right]
-\operatorname{sym}\!\left[(\Hlog\bm{Y}_f)\,\bm{d}_f^{\mathsf T}\right].
\end{equation}
Continuously the second term vanishes because $\bm{H}\bm{Y}=0$. That identity
holds symbolically for every $N$ tested, and it has a measurable consequence:
over $20\,000$ random states per species count, the Fickian part alone and the
full matrix fail the criterion on exactly the same samples, $3$, $13$ and $24$
of $20\,000$ at $N=3,4,5$. The continuous loss of definiteness therefore belongs
to the Fickian part. Discretely, exact conservation for arbitrary
$\bm{d}_f$ and $\jump{\bm{Y}}$ requires
$\bm{1}^{\mathsf T}\bm{Y}_f=1$. Entropy neutrality must be posed on the same
tangent space as the face certificate: the correction contribution vanishes
for every $\bm{z}\in T$ and every $\bm{d}_f$ if and only if
$\bm{Q}^{\mathsf T}\Hlog\bm{Y}_f=0$.

\begin{proposition}[Kernel of the log-mean Hessian]
\label{prop:kernel}
With $\zeta:=\bm{b}^{\mathsf T}\Ylog/Z^{\log}$ and
$s_{\log}:=\bm{1}^{\mathsf T}\Ylog$,
\begin{equation}
\Hlog\Ylog=R_u(1-\zeta)\,\bm{b},
\qquad
\ker\Hlog=
\begin{cases}
\operatorname{span}\{\Ylog\}, & \zeta=1,\\
\{0\}, & \zeta\neq1 .
\end{cases}
\end{equation}
\end{proposition}

\begin{proposition}[Conservative tangent-neutral weight]
\label{prop:tangentweight}
For the closure \eqref{eq:closure}, conservation and entropy neutrality of the
correction term for every admissible jump determine one face weight. When
$\zeta<1$, let $\chi=Z^{\log}(1-\zeta)/s_{\log}$. Its components are
\begin{equation}
Y_{f,s}^{\star}=
\frac{Y_s^{\log}(1+M_s\chi)}
{s_{\log}+\chi\sum_k M_kY_k^{\log}}.
\label{eq:tangentweight}
\end{equation}
Then $\bm{Y}_f^\star>0$, $\bm{1}^{\mathsf T}\bm{Y}_f^\star=1$, and
$\Hlog\bm{Y}_f^\star=c_f\bm{1}$ for a non-negative scalar $c_f$, positive
when $\zeta<1$. At $\zeta=1$ the continuous extension is
$\bm{Y}_f^\star=\Ylog/s_{\log}$ and $c_f=0$.
The uniqueness statement is specific to the rank-one structure of the
ideal-mixture Gibbs Hessian used here and does not imply uniqueness for
arbitrary multicomponent transport closures.
\end{proposition}

Strict concavity and homogeneity of the logarithmic mean give $0<\zeta\le1$,
with equality only for identical positive simplex states. Equation
\eqref{eq:tangentweight} is the regular form of
$(\Hlog)^{-1}\bm{1}/(\bm{1}^{\mathsf T}(\Hlog)^{-1}\bm{1})$. It is therefore
unique when $\zeta<1$; at $\zeta=1$ the kernel in
Proposition~\ref{prop:kernel} gives the same uniqueness after normalization.
For any $\bm{z}\in T$, the correction contribution is proportional to
$(\bm{z}^{\mathsf T}\Hlog\bm{Y}_f^\star)
(\bm{d}_f^{\mathsf T}\bm{z})=0$.

\begin{corollary}[Fickian reduction on the tangent space]
\label{cor:tangent_reduction}
For the weight $\bm{Y}_f^\star$ of Proposition~\ref{prop:tangentweight}, define
\begin{equation}
\bm{K}^{\mathrm F}_f=\operatorname{sym}\!\left[\Hlog\operatorname{diag}(\bm d_f)\right].
\end{equation}
Then
\begin{equation}
\bm Q^{\mathsf T}\bm K^h_f\bm Q
=\bm Q^{\mathsf T}\bm K^{\mathrm F}_f\bm Q,
\qquad
\lambda_{\star,f}^h(\bm d_f;\bm Y_f^\star)
=\lambda_{\star,\mathrm F,f}^h(\bm d_f),
\label{eq:tangent_reduction}
\end{equation}
where the second eigenvalue uses the same denominator
$\bm Q^{\mathsf T}\Hlog\bm Q$. Thus the tangent-neutral correction weight neither
improves nor worsens the Fickian face criterion; it removes only the rank-one
correction contribution on admissible jumps.
\end{corollary}

The distinction from an ambient-space kernel condition is substantive. For
$\bm{Y}_L=(0.2,0.3,0.5)$, $\bm{Y}_R=(0.4,0.1,0.5)$ and
$\bm{M}=(2,16,28)$, one obtains $\zeta=0.9815703\ne1$ and
$\bm{Y}_f^\star=(0.2822101,0.1863656,0.5314243)$, while
$\Hlog\bm{Y}_f^\star=0.0271125\bm{1}$. Thus the correction is neutral on the
physical tangent space although $\Hlog\bm{Y}_f^\star\ne0$.

Other conservative weights need not be tangent neutral. For example, the
normalized log-mean $\bm{Y}_f=\Ylog/s_{\log}$ has the following correction
contribution to entropy production. Here $\Delta\sigma_f$ denotes the
entropy-production contribution, namely
the negative of this face's contribution to $\mathrm d\Phi/\mathrm dt$.
\begin{equation}
\Delta\sigma_f=-\frac{\rho_f|f|}{d_{n,f}}\cdot\frac{R_u(1-\zeta)}{s_{\log}}
\left(\bm{d}_f^{\mathsf T}\jump{\bm{Y}}\right)\!
\left(\bm{b}^{\mathsf T}\jump{\bm{Y}}\right).
\label{eq:defect}
\end{equation}
The closed form was checked on $2\,000$ random states to a maximum relative
error of $8.1\times10^{-11}$. Its sign is indefinite. Since
$1-\zeta=O(\|\jump{\bm{Y}}\|^2)$, for a smooth field with
$\jump{\bm{Y}}=O(h)$, this particular defect is $O(h^2)$ relative to an $O(1)$
leading term. It is a property of the chosen weight rather than an
unavoidable incompatibility.

\subsection{Face weighting and when the criterion is activated}

\begin{figure}[htbp]
  \centering
  \begin{subfigure}{0.48\textwidth}\caption{}%
    \includegraphics[width=\textwidth,height=0.19\textheight]{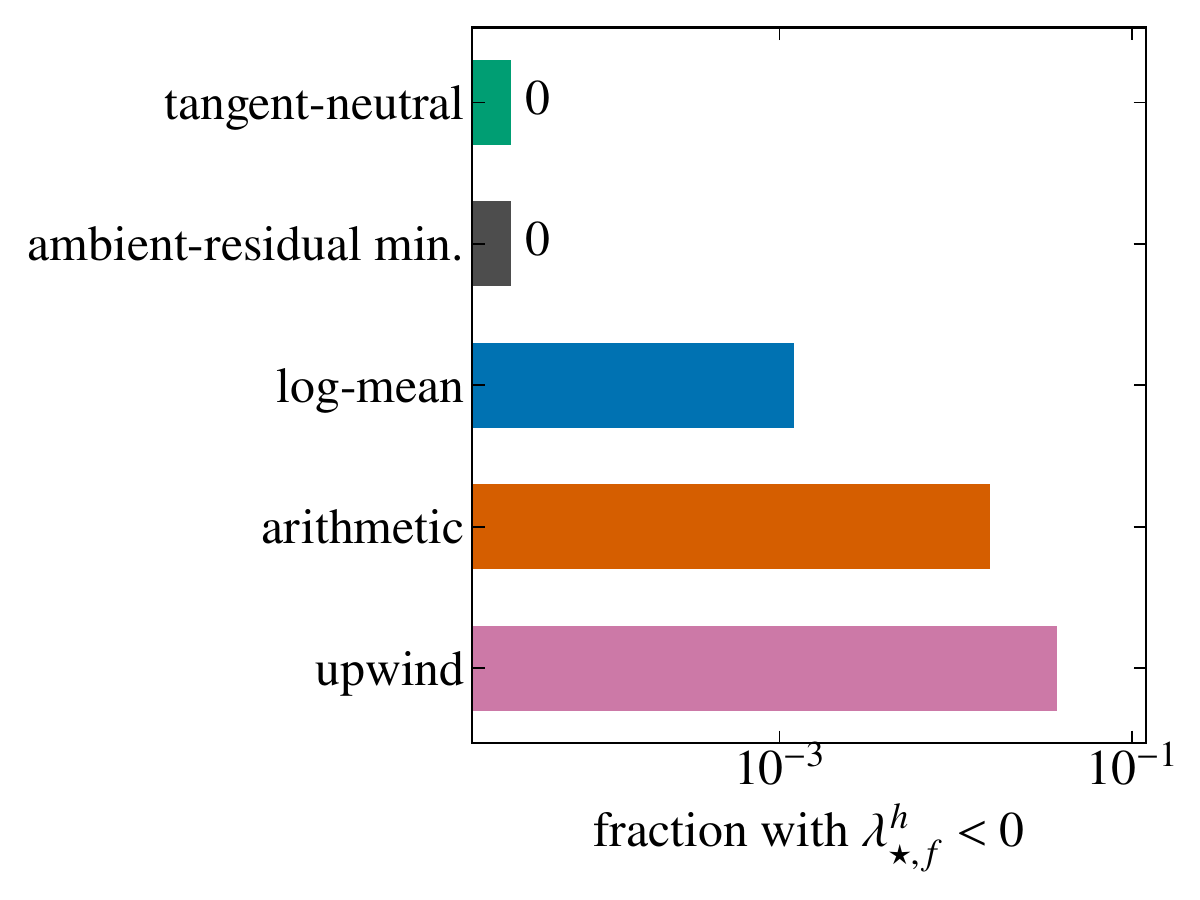}%
    % [图源] figures/scripts/fig1_face_closure.py | 数据: data/p59_tangent_neutral_weight_population.csv
    \label{fig:face_a}\end{subfigure}\hfill
  \begin{subfigure}{0.48\textwidth}\caption{}%
    \includegraphics[width=\textwidth,height=0.19\textheight]{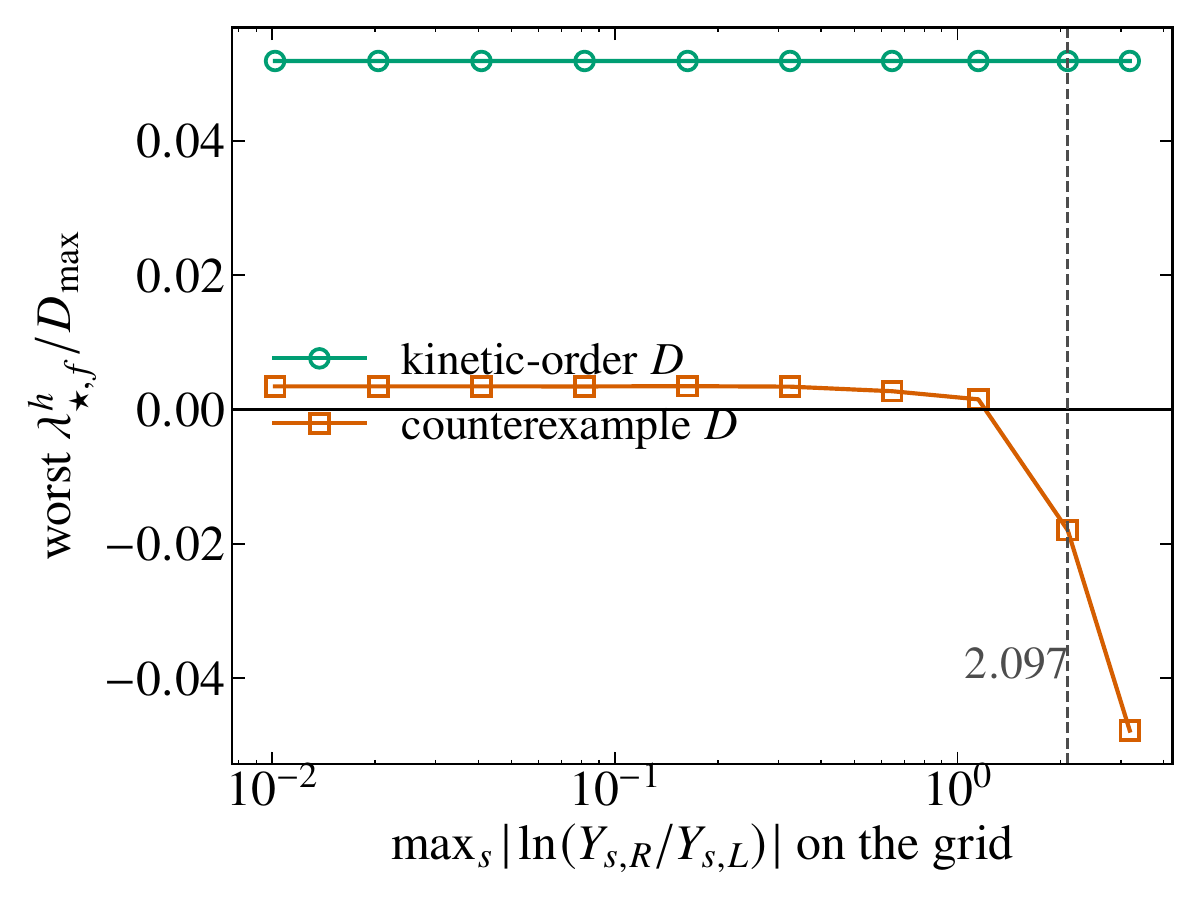}%
    % [图源] figures/scripts/fig1_face_closure.py | 数据: data/a3_resolution.csv
    \label{fig:face_b}\end{subfigure}
  \caption{The two-point criterion and its activation.
  (\textit{a}) Fraction of $39\,999$ random face states violating
  $\lambda_{\star,f}^h\ge0$ for five face weights; the tangent-neutral and
  ambient-residual-minimising weights have no failures and are plotted at the
  axis floor. All five weights conserve mass to
  better than $5\times10^{-15}$.
  (\textit{b}) Worst generalised margin on each refinement level against the
  largest per-species log jump carried by that level, for the counterexample
  diffusivity vector and for a kinetic-order transport model. Here
  $D_{max}:=\max_s d_{s,f}$, so the plotted ordinate is
  $\lambda_{\star,f}^h/D_{max}$. The dashed line marks the critical jump
  $2.097$. The comparison separates activation of the algebraic criterion
  from failure caused by the choice of face weight.}
  \label{fig:face}
\end{figure}

Figure~\ref{fig:face_a} reports the failure fraction over $39\,999$ random face
states for a ternary mixture: $1.558\times10^{-2}$ for the arithmetic mean,
$1.200\times10^{-3}$ for the normalised log-mean, and none for either the
tangent-neutral weight or the historical weight that minimises the stronger
ambient residual $\|\Hlog\bm{Y}_f\|$. The tangent-neutral weight has a worst
normalised margin of $0.01527$, compared with $0.01273$ for the ambient
minimizer. Its smallest component is $9.63\times10^{-7}$ at a state floor of
$10^{-6}$, and the mass residual is below $5\times10^{-15}$ for all five
weights. Conservation alone therefore does not distinguish the weights; the
tangent-space condition does.

Corollary~\ref{cor:tangent_reduction} fixes the interpretation of this result:
for the tangent-neutral weight, the full face margin is exactly the Fickian
margin evaluated with $\Hlog$ at the log-mean state. The zero failure fraction
therefore applies to this ternary mixture and diffusivity population rather
than to arbitrary diffusivities. On the same
$20\,000$ samples per species count used in the independent symbolic check,
the Fickian and tangent-neutral margins are negative on exactly the same
$3$, $13$ and $24$ faces for $N=3,4,5$, respectively. These failures therefore
belong to the Fickian criterion and are independent of the tangent-neutral
weight. The same frozen ternary population gives a $3.72\%$ failure fraction
for the upwind weight in Figure~\ref{fig:face_a}.

\begin{proposition}[Exact shift]
\label{prop:shift}
For any weight with $\bm{1}^{\mathsf T}\bm{Y}_f=1$ and any $\delta\ge0$,
$\lambda_{\star,f}^h(\bm{d}+\delta\bm{1})=\lambda_{\star,f}^h(\bm{d})+\delta$.
\end{proposition}

Consequently $\delta_f=\max(0,-\lambda_{\star,f}^h)$ is the smallest common
diffusivity that restores the face inequality while staying inside the same
conservative closure family, and it vanishes identically on faces that already
pass. The repair is a minimal common artificial diffusivity, equivalent here to
adding a Fickian term $\delta\rho\nabla\bm Y$; the measured shift error is below
$10^{-16}$.

The random diffusivity vectors behind Figure~\ref{fig:face_a} are synthetic and
are not generated by the transport controls considered below. With the
arithmetic face composition, the discrete criterion is satisfied on all
$5\,210$ faces extracted from production neutral-air/carbon profiles, with a
smallest normalised margin of $+0.4888$, and on all $36\,000$ samples in the
quasi-neutral ambipolar-air control, with a smallest margin of $+0.161$.
Refining a composition front gives a complementary comparison
(Figure~\ref{fig:face_b}). The criterion is violated
only once a face carries a per-species log jump above $2.097$, a mass-fraction
change of about a factor of eight, whereas the kinetic-order model holds a
worst margin of $+0.05187$ on every level and survives jumps of
$|\Delta\ln Y|=18.4$.

The criterion of Corollary~\ref{cor:telescope} is therefore exact. General face
weights can change it, whereas the tangent-neutral weight leaves the Fickian
margin unchanged by Corollary~\ref{cor:tangent_reduction}. In the populations
examined here, the weight dependence is activated by the synthetic diffusivity
vectors but not by either physically constrained transport control.
Sections~\ref{sec:layer2} and
\ref{sec:layer3} concern layers where the corresponding failures are activated.

\section{Layer 2: multipoint non-orthogonal reconstruction}
\label{sec:layer2}

\subsection{The per-face certificate does not transfer}

The face certificate of Section~\ref{sec:layer1} is a statement about a
two-point stencil. A tangential correction makes the face flux depend on a
wider stencil, and the quadratic form behind
Corollary~\ref{cor:telescope} no longer decomposes face by face.

The WLS prototype behaves as intended on the properties it was built for. On
clipped Voronoi polygons it reproduces constant and linear fields to
$2.0\times10^{-15}$ and reduces the uncorrected normal-gradient error from
$0.937$ to $3.5\times10^{-15}$. Conservation is unaffected: the face mass
residual stays at $1.3\times10^{-15}$ and the global species residual at
$2.1\times10^{-19}$, at a worst WLS condition number of $2.83$. Over
$55\,400$ face-state samples, however, $540$ local entropy contributions are
positive. The per-face
certificate is therefore lost, and global dissipation over $400$ sampled state
fields is an observation about those fields rather than a certificate.

\begin{figure}[htbp]
  \centering
  \begin{subfigure}{0.48\textwidth}\caption{}%
    \includegraphics[width=\textwidth,height=0.19\textheight]{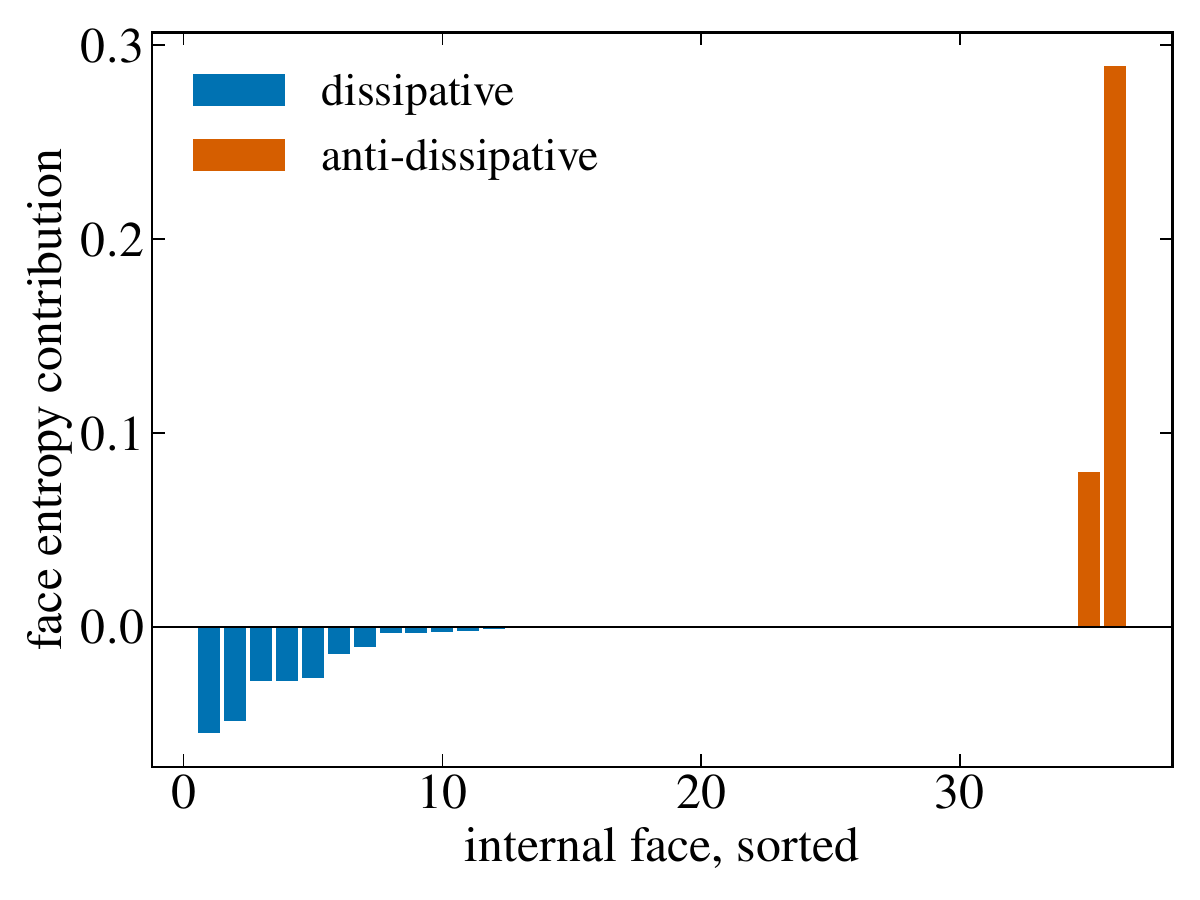}%
    % [图源] figures/scripts/fig2_wls_entropy.py | 数据: data/p35_wls_global_entropy_counterexample.json
    \label{fig:wls_a}\end{subfigure}\hfill
  \begin{subfigure}{0.48\textwidth}\caption{}%
    \includegraphics[width=\textwidth,height=0.19\textheight]{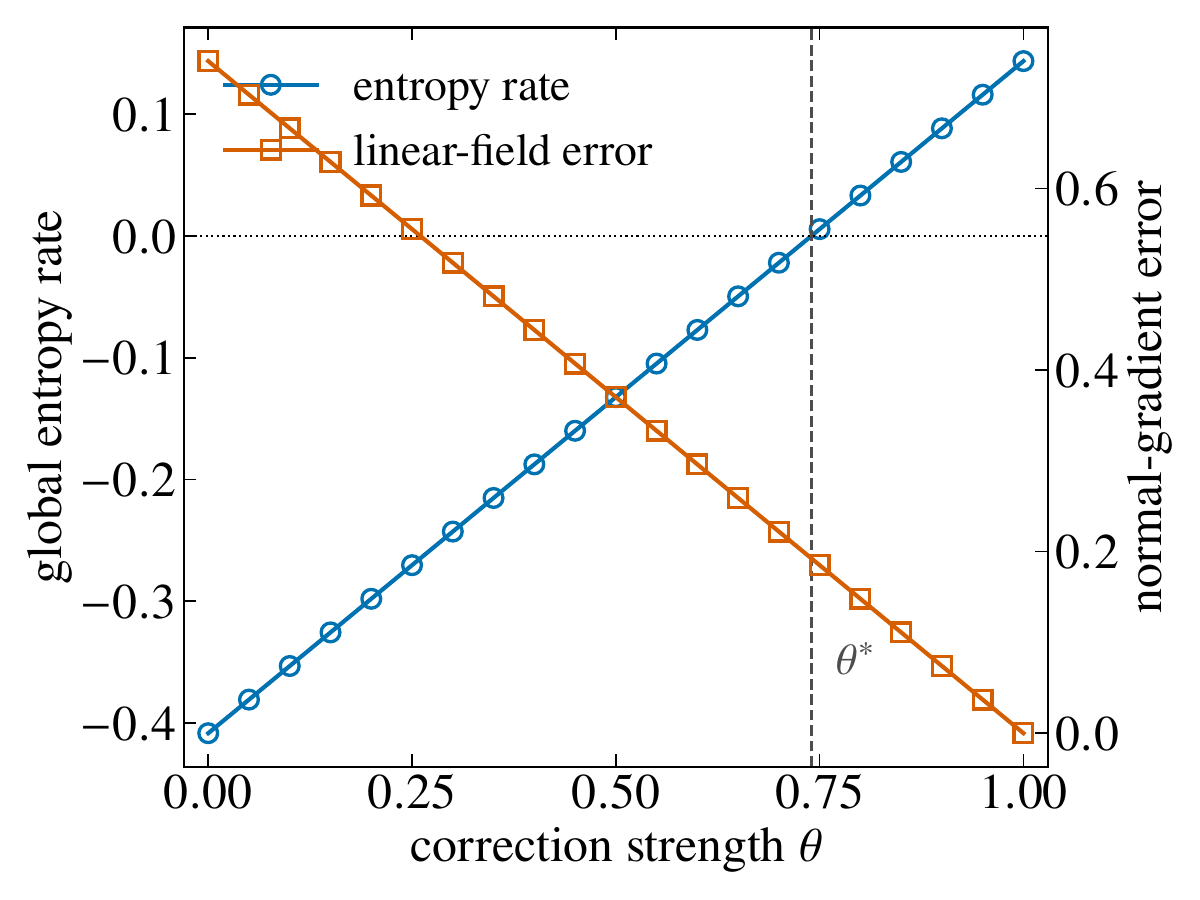}%
    % [图源] figures/scripts/fig2_wls_entropy.py | 数据: data/p38_wls_theta_tradeoff.json
    \label{fig:wls_b}\end{subfigure}
  \caption{Local and global entropy behaviour of the WLS correction.
  (\textit{a}) Sorted face entropy contributions at the counterexample state on
  a $16$-cell mesh; five of $36$ internal faces are anti-dissipative and
  dominate the sum.
  (\textit{b}) Global entropy rate and linear-field normal-gradient error
  against the correction strength $\theta$. The rate is affine in $\theta$ and
  crosses zero at $\theta^\ast=0.7398$, where the linear-field error is
  $0.193$.}
  \label{fig:wls}
\end{figure}

Losing the local certificate does not by itself imply that the global rate can
turn positive. To decide that, the diffusivities were made equal, which renders
the projected two-point member unconditionally dissipative and removes the
Layer~1 mechanism from the experiment. A search over $600$ states and four
optimised starts returns a state at which the global entropy rate is
$+0.1436$, against $-0.4083$ for the projected two-point member defined by
\eqref{eq:reconstruction} at $\theta=0$, with an entropy-telescoping residual
of $1.1\times10^{-16}$. All $32$ perturbations with logit standard deviation
$10^{-3}$ remain above the counterexample threshold, showing local persistence
in the chosen parameterisation.
Figure~\ref{fig:wls_a} shows that five of the $36$ internal faces carry the
positive contribution.

As an implementation-independence check, the same frozen state and geometry
were reassembled by a NumPy-only script that imports none of the project
geometry, thermodynamics, diffusion or WLS modules. It regenerates the clipped
Voronoi cells from the fixed archived seed and reads only the archived state. The
36-face assembly reproduces $-0.4083058819$ for $\theta=0$ and
$+0.1436029942$ for $\theta=1$ (both differences from the reference are zero
at double precision); the largest face mass residual is
$4.1\times10^{-16}$ and the largest global species residual is
$5.9\times10^{-21}$. This independently reproduces the reported 16-cell
counterexample.

To test sensitivity to the distance weighting, the reconstruction was rebuilt
for weights
$\|\bm{r}_{ij}\|^{-\kappa}$ with $\kappa=0,1,2,3$ and repeating the search on the same
geometry answers the narrower question. Every weighting is linearly exact to
$6.2\times10^{-16}$ and conservative, and every one admits a positive global
rate, between $+0.126$ and $+0.163$, with $5$ to $10$ anti-dissipative faces.
The inverse-square case reproduces the value quoted above, so the counterexample
is not confined to $\kappa=2$ within this tested family.

An independent search on $32$, $64$ and $96$ cell meshes returned no positive
rate (best normalised rates $-0.845$ to $-0.938$), a negative finite-search
result rather than evidence of mesh-independent stability.

\subsection{Accuracy against entropy}

The counterexample state also exposes the trade-off directly. The global
entropy rate is affine in $\theta$ to a residual of $1.1\times10^{-16}$, rising
from $-0.408$ at $\theta=0$ to $+0.1436$ at $\theta=1$ and crossing zero at
$\theta^\ast=0.7398$. The linear-field normal-gradient error runs the other
way, from $0.740$ at $\theta=0$ to machine zero at $\theta=1$, and is still
$0.193$ at $\theta^\ast$ (Figure~\ref{fig:wls_b}). For this prototype and this
state, the correction strength that makes the operator linearly exact is
exactly the one that makes it anti-dissipative.

The affine structure suggests treating $\theta$ as a global, state-dependent
control. Capping the entropy rate alone is not enough. Applied to an explicit
Euler trajectory from the counterexample state, the cap does keep positivity,
conservation and potential decay on every accepted step. The right-hand side
nevertheless still points outward at the simplex boundary. After $80$ halvings
the step collapses to $8.3\times10^{-27}$, having covered $0.99\%$ of the target
interval. The entropy-only control therefore does not provide a viable positive
trajectory in this test. Intersecting the entropy constraint with
the one-step positivity constraint gives non-empty joint intervals
$[0,0.4187)$ and $[0,0.4158)$ at the two states examined, both bounded by a
single species in a single cell, and both excluding the entropy-only cap.
Solving for a common $\theta$ at every step then advances $100$ steps of
$\Delta t=10^{-4}$ with no rejections and a total potential decrease of
$1.77\times10^{-3}$. The accuracy cost is explicit: the manufactured linear-gradient
error along that trajectory is $0.219$ to $0.432$. The construction is a
global, state-dependent diagnostic, not a local limiter and not an
invariant-domain result.

\subsection{Consistency depends on mesh geometry}
\label{sec:layer2:geom}

Accuracy of the corrected operator was assessed on manufactured solutions.
A strictly interior affine steady field is reproduced to machine zero on four
meshes and three time levels with $\theta=1$ throughout. A curved cosine field
is not: the largest spatial residual is $2.7358\times10^{-4}$, exceeding the
reference tolerance of $5\times10^{-5}$; the corresponding accuracy criterion
is therefore not met.
Paired-grid decomposition
attributes part of the residual to the difference between centroid point values
and polygonal cell averages, which reaches $3.35\times10^{-3}$. On regular
orthogonal meshes with analytic cell averages the same manufactured field
converges at observed orders not below $1.986$ in $L_2$ and $1.965$ in
$L_\infty$, which narrows the failure to the combination of the tested WLS
prototype with random Voronoi geometry. The orthogonal control is not used to
alter this conclusion.

\begin{figure}[htbp]
  \centering
  \begin{subfigure}{0.48\textwidth}\caption{}%
    \includegraphics[width=\textwidth,height=0.19\textheight]{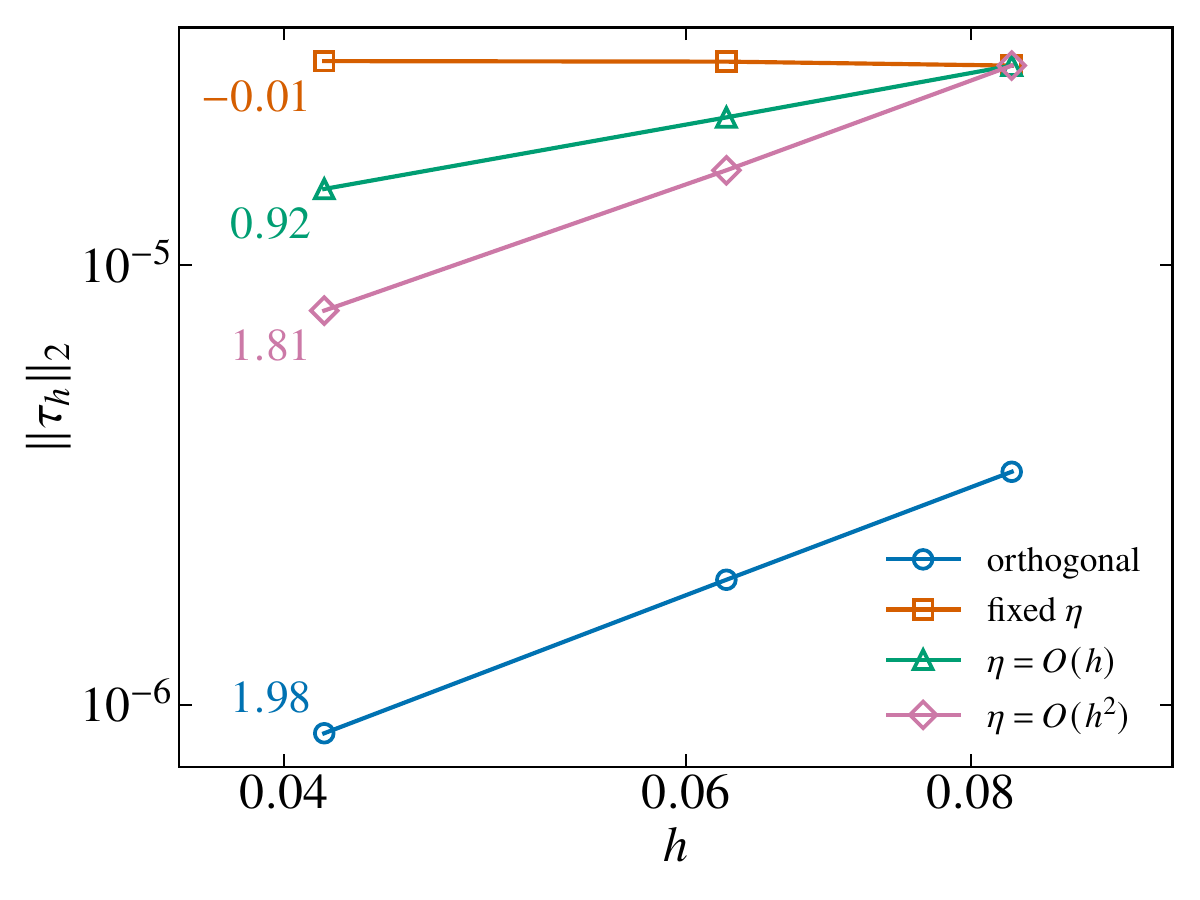}%
    % [图源] figures/scripts/fig3_geometric_scaling.py | 数据: data/p48_nonorthogonality_scaling_audit.json
    \label{fig:scaling_a}\end{subfigure}\hfill
  \begin{subfigure}{0.48\textwidth}\caption{}%
    \includegraphics[width=\textwidth,height=0.19\textheight]{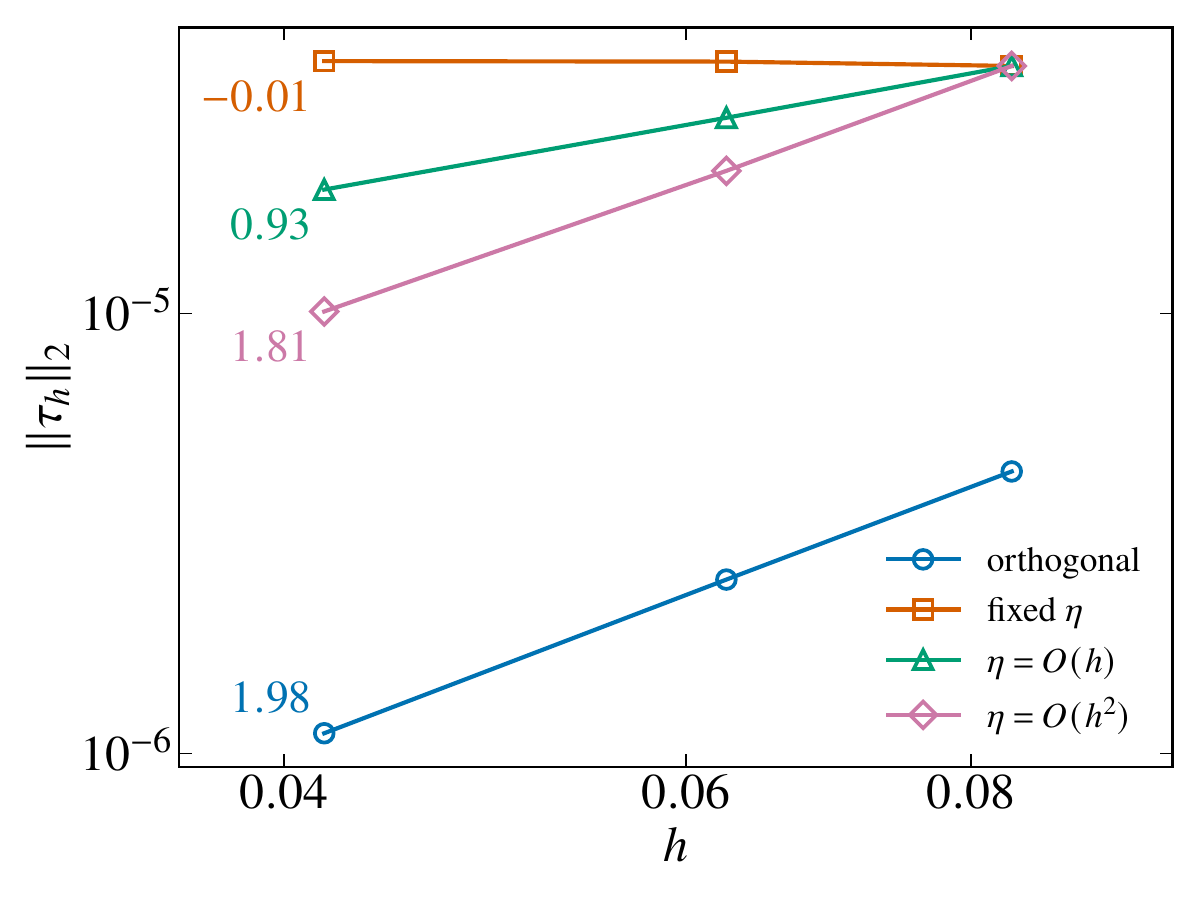}%
    % [图源] figures/scripts/fig3_geometric_scaling.py | 数据: data/p50_state_dependent_multicomponent_mms.json
    \label{fig:scaling_b}\end{subfigure}
  \caption{Truncation error against mesh spacing for four non-orthogonality
  families, labelled with the observed $L_2$ order on the finest interval.
  (\textit{a}) Common diffusivity.
  (\textit{b}) Unequal, state-dependent diffusion with a
  largest-to-smallest diffusivity ratio of $19.95$ and a mass-correction flux
  reaching $58.0\%$ of the Fickian flux. The classification is unchanged.}
  \label{fig:scaling}
\end{figure}

A controlled study separates the geometric mechanism. Displacing a regular
lattice by a fixed fraction of $h$ amplifies the residual by $15.1$ in $L_2$
and $25.0$ in $L_\infty$ on the finest mesh, with non-orthogonality correlating
with the $L_2$ residual above $0.988$ on every mesh. Its corresponding diagnostic
has an unperturbed $L_\infty$ order of $1.818$, below the reference value
$1.9$. Scaling the displacement then produces the
ordering shown in Figure~\ref{fig:scaling_a}: observed $L_2$ orders of $1.98$
for the orthogonal control, $-0.01$ for fixed relative displacement, $0.92$ for
$O(h)$ and $1.81$ for $O(h^2)$.

\paragraph{Conditional scaling argument}
A quadratic-field expansion suggests that, under paired quasi-uniform
perturbation of a second-order orthogonal family, uniformly invertible WLS
moment matrices with uniformly bounded stencils, a smooth field with bounded
derivatives, and a consistency-compatible boundary closure,
where $\bm{\tau}_h$ denotes the truncation residual of the corrected spatial
operator,
\begin{equation}
\|\bm{\tau}_h\|=O(h^2+\eta_h) .
\label{eq:trunc}
\end{equation}
The local expansion and face contraction identity are verified to residuals
below $4.4\times10^{-15}$, and the predicted exponents differ from the measured
ones by at most $0.256$. Equation~\eqref{eq:trunc} remains a working scaling
hypothesis, however: the norm defining $\bm{\tau}_h$, boundary contribution,
and $h$-uniform constants $C_0,C_1$ have not been established. It classifies
the observed orders but is neither a quantitative bound nor an entropy
statement.

Repeating the study with unequal, state-dependent diffusion leaves the
classification intact (Figure~\ref{fig:scaling_b}). At a
largest-to-smallest diffusivity ratio of $19.95$, with the mass-correction flux
reaching $58.0\%$ of the Fickian flux, the four families give finest $L_2$
orders of $1.983$, $-0.006$, $0.926$ and $1.814$. The state-dependent matrix
introduces no additional $O(h)$ term at these settings.

On the same twelve states, the global entropy rate of the full correction is
negative throughout, between $-5.52\times10^{-4}$ and $-5.38\times10^{-4}$, and
both the global cap and the dynamic common $\theta$ stay inactive at
$\theta=1$. The non-orthogonal families nevertheless produce $166$ positive
local face contributions in total, up to $48$ in a single state. Global
dissipation on smooth interior states does not restore the per-face
certificate.

\subsection{Paired refinement and cross-seed sufficiency}
\label{sec:layer2:robust}

Unlike Figure~\ref{fig:scaling}, which classifies geometric scaling,
the paired-refinement control separates temporal from spatial error.
The supporting plots are collected in \ref{app:verification}
(Figure~\ref{fig:pairing}). Three grids
with four Heun time levels each give a smallest observed temporal
self-convergence order of $1.99991$ and a Richardson time-contamination ratio
of at most $2.04\times10^{-6}$ on the finest diagonal, so the diagonal spatial
orders are not polluted by time error. Those orders are $1.983$ and $1.949$ for
the orthogonal family and $2.676$ and $1.857$ for the $O(h^2)$ family
(Figure~\ref{fig:pairing_a}). All $24$ trajectories keep $\theta=1$, conserve
species to $6.3\times10^{-17}$, stay positive at a smallest mass fraction of
$0.0851$, and decrease the potential at every step.
The value $2.676$ is a finite three-grid estimate; the available data do not
identify whether its elevation reflects pre-asymptotic behaviour, cancellation,
or another finite-grid effect. It is retained as an observed order and is not
used as a theorem-level claim about the asymptotic order of the corrected
operator.
The paired space--time refinement separates temporal from spatial error, and
the multi-seed population preserves the ordering
$\text{fixed }\eta>O(h)>O(h^2)>\text{orthogonal}$ across the tested modes and
norms. The supporting plots and complete numerical record are in
\ref{app:verification}.
The multi-seed population does not meet the interval-stability target and is
therefore used as a sensitivity check rather than as evidence of asymptotic
convergence.

\section{Layer 3: positivity-preserving time integration}
\label{sec:layer3}

\subsection{Positivity and conservation do not imply dissipation}

Section~\ref{sec:layer1} concerns the spatial operator alone. Time integration
imposes a separate requirement, and for stiff chemistry the schemes of choice
are those that keep positivity at any step size. In this section
$\bm f(\bm Y)=\bm\omega(\bm Y)/\rho$ denotes the source-only rate.
For the source-only comparisons below, the reported potential is the specific
potential $\varphi$ of one unit-volume cell; the displayed factor $m/R_u$
renders the quoted gaps dimensionless.
For a time step, write $\Delta\bm Y:=\bm Y^{n+1}-\bm Y^n$.
For the single-network step-size comparison below, define
\[
\tau_{\mathrm{chem}}
=\left[\max_{\lambda\in\operatorname{spec}(J_0)}|\lambda|\right]^{-1},
\qquad J_0=\left.\frac{\partial\bm f}{\partial\bm Y}\right|_{\bm Y_0}.
\]
The archived calculation evaluates $J_0$ by forward differences,
with increment $10^{-8}\max(Y_{0,s},10^{-12})$ in coordinate $s$.
This initial-state spectral timescale is held fixed throughout that scan;
it is not recomputed along each trajectory.

\begin{proposition}[Ascent barrier]
\label{prop:ascent}
Let $\Phi$ be convex and differentiable, $\Psi$ any one-step map, and
$\Delta\bm y:=\Psi(\bm{y})-\bm{y}$. If
$\nabla\Phi(\bm{y})\cdot\Delta\bm y\ge0$ then
$\Phi(\bm{y}+\gamma\Delta\bm y)\ge\Phi(\bm{y})$ for every $\gamma\ge0$ such that
$\bm{y}+\gamma\Delta\bm y\in\operatorname{dom}\Phi$.
\end{proposition}

Convexity gives the result immediately. Its consequence is structural: when the
update direction ascends, every positive scalar relaxation along that same
fixed direction leaves the increase in place while the relaxed state remains
in the domain. The obstruction is one of direction, not of the scalar distance
travelled along that ray.

\begin{figure}[htbp]
  \centering
  \begin{subfigure}{0.48\textwidth}\caption{}%
    \includegraphics[width=\textwidth,height=0.19\textheight]{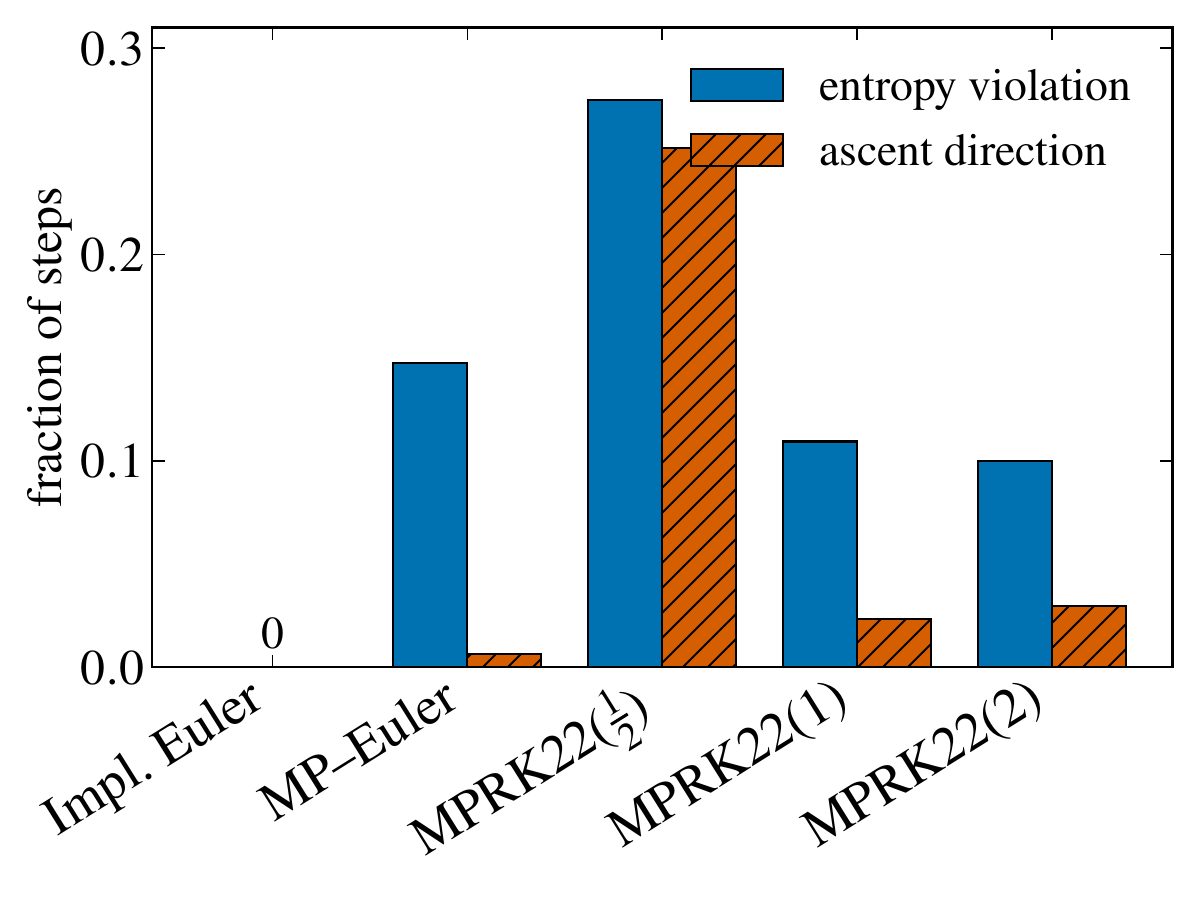}%
    % [图源] figures/scripts/fig5_patankar_barrier.py | 数据: data/p0_barrier.csv; data/p61_implicit_euler_diagnostic.json
    \label{fig:barrier_a}\end{subfigure}\hfill
  \begin{subfigure}{0.48\textwidth}\caption{}%
    \includegraphics[width=\textwidth,height=0.19\textheight]{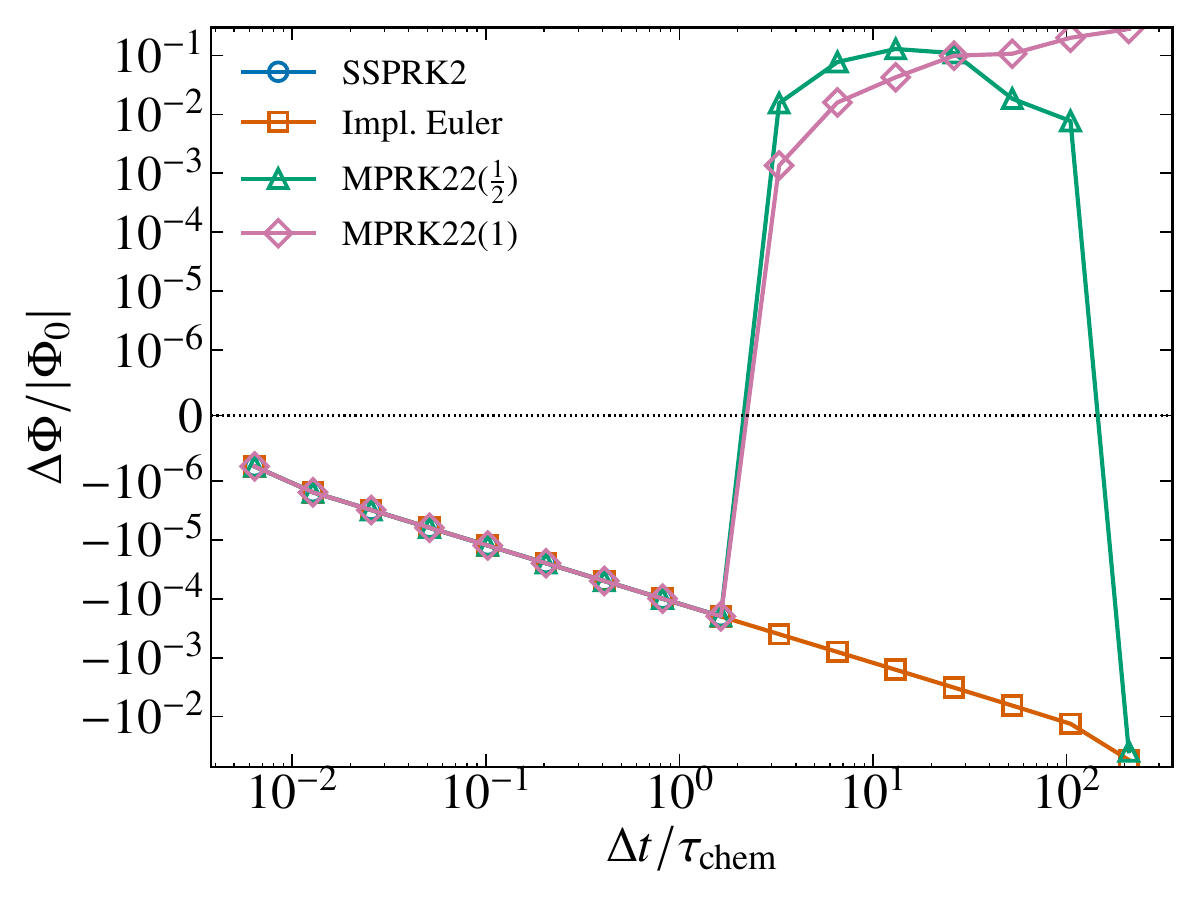}%
    % [图源] figures/scripts/fig5_patankar_barrier.py | 数据: data/a2_trilemma.csv
    \label{fig:barrier_b}\end{subfigure}
  \caption{The Patankar entropy barrier.
  (\textit{a}) Fraction of steps that raise the free energy, and fraction that
  do so along an ascent direction, over $150$ networks and $30$ step sizes.
  Among the $4\,408$ converged positive implicit-Euler roots, none has either
  property. First-order MP--Euler violates on $14.76\%$ of steps.
  (\textit{b}) Worst relative potential change against step size on a single
  network, restricted to steps that stay positive. Classical schemes lose
  positivity before they raise $\Phi$; the Patankar schemes stay positive
  throughout and raise $\Phi$ above $\Delta t\approx3\tau_{\mathrm{chem}}$.}
  \label{fig:barrier}
\end{figure}

MP--Euler is the first-order Patankar update used here. SSPRK2 and SSPRK3 are
strong stability-preserving Runge--Kutta controls of orders two and three.
MPRK22 denotes the second-order modified Patankar--Runge--Kutta family.
Figure~\ref{fig:barrier_a} shows that the discriminating feature is the
Patankar rescaling rather than the order. Over $150$ networks and $30$ step
sizes spanning five decades, $4\,408$ of $4\,500$ implicit-Euler solves meet
the nonlinear residual tolerance and have positive roots; none raises the
free energy or follows an ascent direction. The remaining $92$ attempts are
solver failures (56 iteration limits and 36 positivity-line-search failures),
not admissible time steps. In contrast, the first-order MP--Euler scheme violates on
$14.76\%$ of all $4\,500$ steps. The second-order MPRK22 family violates on
$9.98$--$27.49\%$ of steps, with ascent-direction fractions of $0.62$--$25.18\%$
across the Patankar family and zero for the classical family. Where a positive
scalar rescaling $\gamma$ can repair a violation, its mean value is $0.54$ to
$0.70$; this changes the proposed update by $30$--$46\%$ and removes the
unconditional character that motivates
the schemes.

Figure~\ref{fig:barrier_b} shows the same on a single network. SSPRK2 and
SSPRK3 lose positivity at $\Delta t/\tau_{\mathrm{chem}}=1.64$ and never raise
$\Phi$ while positive. The converged implicit-Euler branch is positive and
dissipative across the sampled range to $209.6$. MPRK22 is positive across the whole range and first
raises $\Phi$ at $\Delta t/\tau_{\mathrm{chem}}=2.96$, violating on $39.6\%$ and
$43.8\%$ of positive runs for $\alpha=1/2$ and $\alpha=1$, with worst relative
increases of $0.134$ and $0.277$ at mass errors below $3\times10^{-14}$. The
violating set is not an interval in $\Delta t$: a single very large step lands
close to equilibrium and does not violate. A critical step size therefore
cannot be found by bisection.

\subsection{An explicit counterexample}
\label{sec:layer3:counterexample}

Ensembles establish frequency, not impossibility. A three-species network
settles the question in closed form. Take
\begin{equation}
2B\rightleftharpoons A+2C,\qquad 2B\rightleftharpoons 2A+C,
\end{equation}
with $\bm{M}=(m,\tfrac32 m,m)$, rate constants $(1,100)$, zero standard
potentials and $\bm{Y}=(6,3,2)/11$. The mass-transfer graph is irreducible.
The infinite-step limit of the MP--Euler update is the normalised kernel vector
of the update matrix, which the matrix-tree formula gives exactly as
$\bm{v}=(3,600,4)/607$; the computed vector matches to $2.3\times10^{-16}$.
Bounding the resulting logarithmic combination by rational
$\operatorname{artanh}$ series gives
\begin{equation}
0.7977206390010004<\frac{m}{R_u}\left[\varphi(\bm{v})-\varphi(\bm{Y})\right]
<0.7977206390010005 .
\end{equation}
The same construction at the finite step $\Delta t=100$ has exact solution
$\bm{Y}^{n+1}=(48102,1719525,20126)/1787753$, which is strictly positive and
exactly conservative, and
\begin{equation}
0.6937764562608661<\frac{m}{R_u}\left[\varphi(\bm{Y}^{n+1})-\varphi(\bm{Y})\right]
<0.6937764562608662 .
\end{equation}
A first-order, unconditionally positive and conservative Patankar scheme
therefore raises the free energy at a finite step size, with a certificate that
does not depend on floating-point arithmetic. Repeating the construction in
rational arithmetic through all stages of MPRK22 gives strictly positive lower
bounds on both $\Delta\Phi$ and $\nabla\Phi^n\cdot\Delta\bm{Y}$ at
$\alpha=1/2$ and $\alpha=1$ with exactly zero mass error, so the increase
occurs along an ascent direction and Proposition~\ref{prop:ascent} applies.
This result does not contradict the positivity guarantees of Patankar methods;
it identifies an additional thermodynamic condition that is not implied by
positivity alone.
The counterexample and the ensemble fractions in this section use the PDS
allocation of Section~\ref{sec:structure}, namely the net reaction rate
$q_r$ followed by the gain/loss projection. They do not cover alternative
splittings that separate forward and reverse reaction rates before the
Patankar rescaling.

\subsection{Coverage of the third-order members}

\begin{figure}[htbp]
  \centering
  \begin{subfigure}{0.48\textwidth}\caption{}%
    \includegraphics[width=\textwidth,height=0.19\textheight]{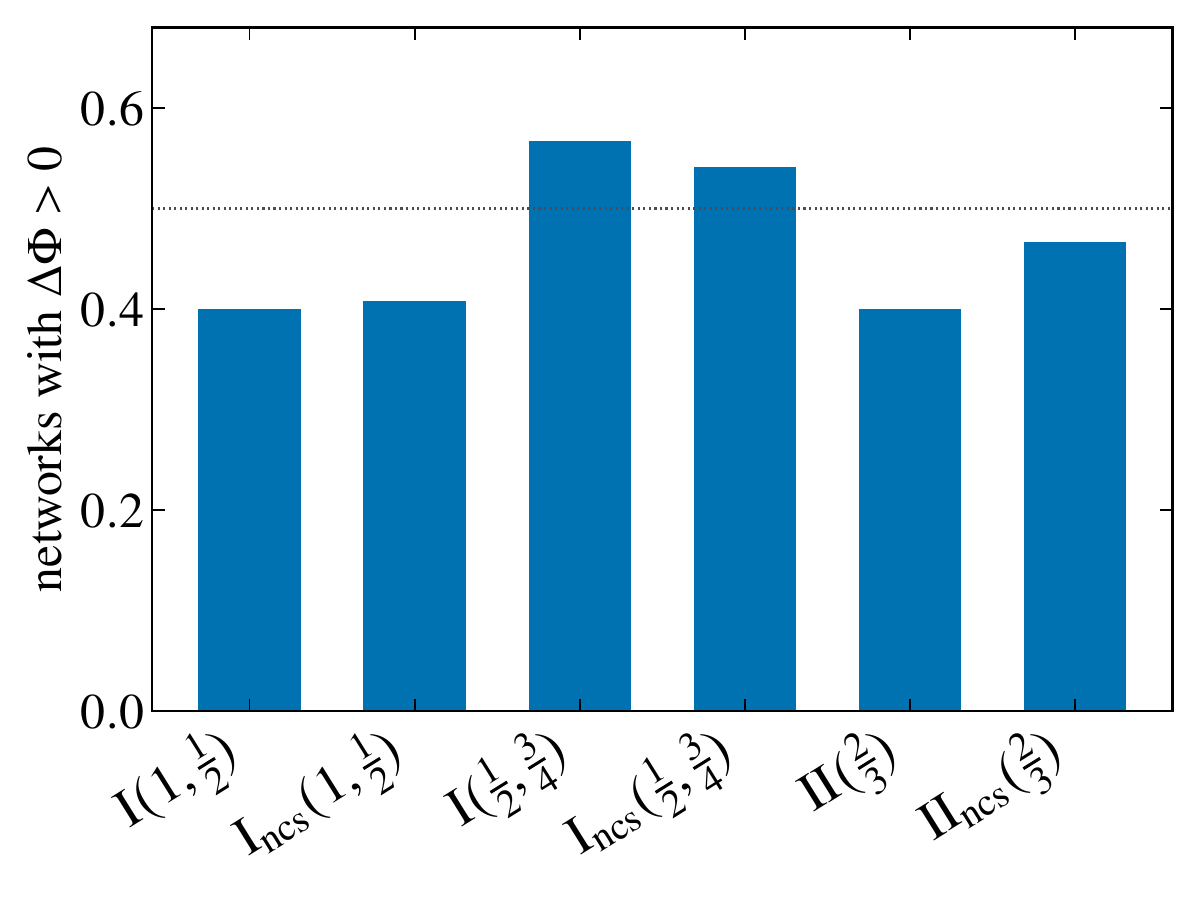}%
    % [图源] figures/scripts/fig6_members_and_relaxation.py | 数据: data/p10_mprk43_ensemble_summary.json
    \label{fig:members_a}\end{subfigure}\hfill
  \begin{subfigure}{0.48\textwidth}\caption{}%
    \includegraphics[width=\textwidth,height=0.19\textheight]{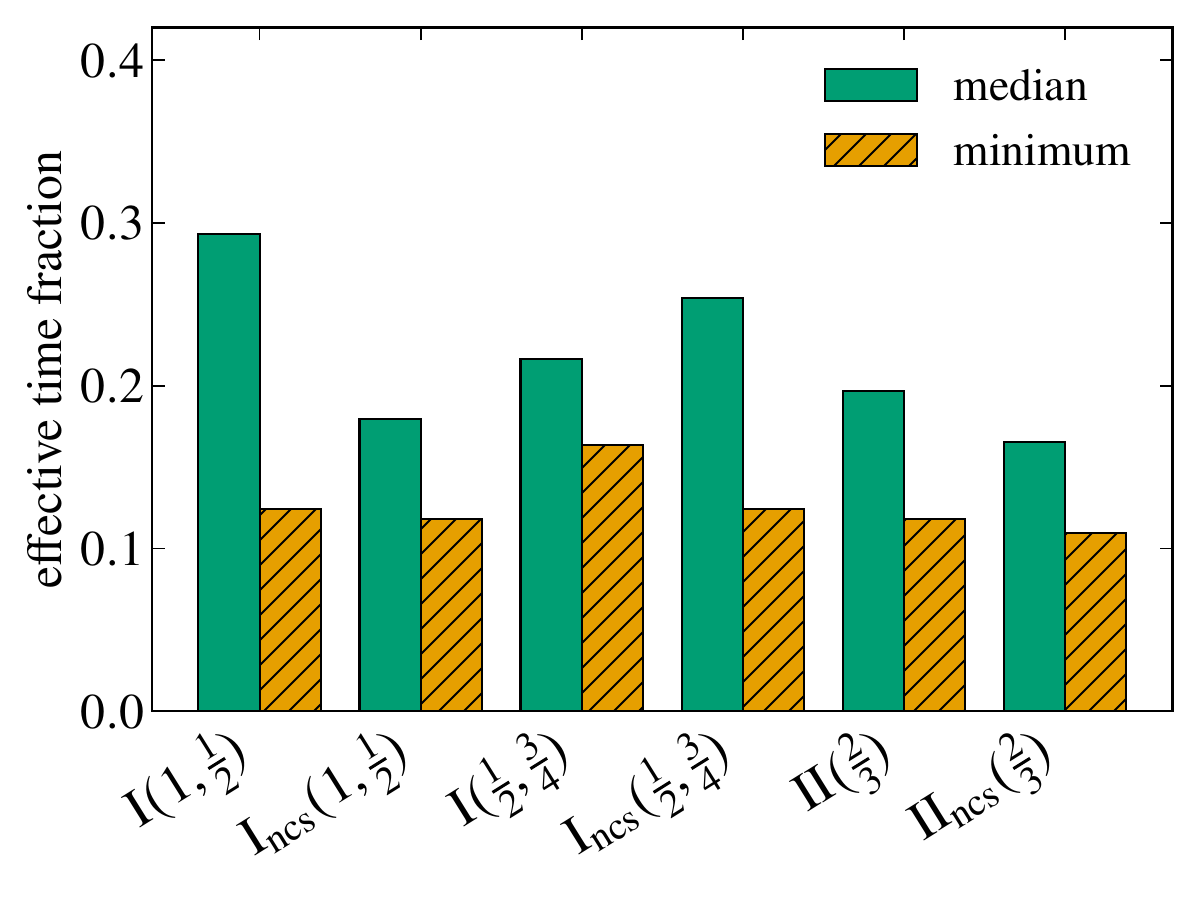}%
    % [图源] figures/scripts/fig6_members_and_relaxation.py | 数据: data/p21_six_member_relaxation_ensemble.json
    \label{fig:members_b}\end{subfigure}
  \caption{Third-order members and the cost of dissipative step control.
  Labels I and II denote MPRK43I and MPRK43II, respectively; their
  parenthesised parameters distinguish the tested members.
The subscript ncs denotes the non-conservative-stage variant
($\delta_{\rm P}=0$); labels without ncs use $\delta_{\rm P}=1$. This distinction concerns internal
  stages, not the conserved final mass.
  (\textit{a}) Fraction of $120$ irreducible detailed-balance networks on which
  each member raises the free energy at some step size. All members remain
  positive and conservative on every step.
  (\textit{b}) Median and smallest fraction of physical time advanced per unit
  of proposed time on $12$ networks, once step rejection and recomputation enforce dissipation.}
  \label{fig:members}
\end{figure}

Here MPDeC denotes modified Patankar deferred correction. The six representative
third-order members of \citet{kopecz2018bit} were
implemented from the source equations and checked against the published
coefficients. All six remain positive and conservative and reach their tested
third-order regime, yet each raises the free energy on a substantial fraction
of the fixed detailed-balance population: $40.0$--$56.7\%$ over $120$ networks
(Figure~\ref{fig:members_a}). An equidistant MPDeC(3) implementation shows the
same qualitative behaviour. These results establish coverage of the named
members, not a theorem for arbitrary Patankar schemes; recent order conditions
for fourth-order MPRK schemes \citep{izgin2024nsark} further widen the family.

\subsection{Step rejection and adaptive recomputation}

Enforcing dissipation by rejecting and recomputing offending steps works, in the sense that trajectories then
satisfy all three constraints, but it converts an unconditional method into an
adaptive one. Applied to a frozen network, the adaptive step-control trajectory
reaches the target time with positivity, conservation and a monotone decrease of
$\Phi$, at a finest observed order of $3.136$, while advancing only $22.7\%$ of
the proposed physical time on the coarse trajectory. Over $24$ networks the
smallest effective time fraction is $0.111$. Extending to all six members over
$12$ independently screened networks
(Figure~\ref{fig:members_b}) gives median effective time fractions of $0.165$ to
$0.293$ and a smallest value of $0.109$, at a largest mass error of
$4.4\times10^{-16}$. Roughly three quarters of the proposed time is discarded
to enforce the inequality.

Order under step rejection and recomputation is not uniform across the ensemble. Six-level refinement
on six networks has errors decreasing strictly and all structural checks passing,
but the finest observed orders for MPRK43I$(1/2,3/4)$ and its non-conservative
stage counterpart are $2.3465$ and $2.7937$, below the reference tolerance
of $2.8$. The value $2.7937$ is therefore treated as a finite-grid
observation rather than an asymptotic order claim; the paired-refinement
classification in Section~\ref{sec:layer2:robust} is intentionally limited,
while the high-precision follow-up
identifies the terminal-floor mechanism for the affected samples. The reference
order criterion is not met. Its mechanism was closed
separately for each sample: a $50$-digit high-precision check recovers order $2.947965$ for one
sample where double precision gives $1.707652$, identifying a terminal error
floor, and confirms $3.0230\to3.0029$ for the other, identifying a
pre-asymptotic transition. A $50$-digit fixed-step RK4 cross-check agrees with
the high-precision reference to $3.7\times10^{-19}$. The diagnosis explains the
failure and does not reverse it.

\paragraph{Independent reproduction}
\label{sec:layer3:thirdparty}

The rational counterexample and the MPRK22 ascent direction were independently
reproduced in the Julia package of \citet{kopecz2025julia}, while a backward
Euler control remained dissipative. The external check assembled the production
matrix directly from the reaction network and called none of the present solver
code; detailed double-precision values are archived in \ref{app:verification}.

\section{What can and cannot be composed}
\label{sec:baseline}

Sections~\ref{sec:layer1}--\ref{sec:layer3} first establish a face-level
certificate and then show two distinct ways in which structural guarantees can
be lost at subsequent discretization layers. What remains is to determine which
of these properties survive fully discrete composition. Each layer supplies a
certificate for a different algebraic object,
and the hypotheses needed to combine them include an invariant-domain
condition that the face entropy inequality does not provide.

The missing implication is exposed by a two-cell finite-volume problem with
one internal face and zero boundary flux. Set $\rho=1$, $V_L=V_R=1$,
$|f|=d_{n,f}=1$, $\bm M=(2,16,28)$ and
$\epsilon=10^{-6}$, $\theta=0$, and the arithmetic face weight. With
$\bm d=(1,5,1)$,
$\bm Y_L=(\epsilon,\epsilon,1-2\epsilon)$ and
$\bm Y_R=(0.1,0.9-\epsilon,\epsilon)$, the arithmetic face weight is conservative
and has $\lambda_{\star,f}^h=0.8562>0$, so no common-diffusivity repair is active.
Nevertheless, the first component of the left-cell semi-discrete right-hand
side is $-0.0800024$. A backward-Euler solve at
$\Delta t=2\times10^{-5}$ gives
$Y_{L,1}^{n+1}=-5.995\times10^{-7}$, with nonlinear residual and total-species
defect both at $1.1\times10^{-16}$. Thus a face entropy certificate does not
make the spatial operator inward pointing at the boundary of the simplex, and
an implicit step does not repair that missing property at this finite step.

The mechanism is a state-dependent step restriction, not failure at every step
size. The step-size study repeats the same solve at prescribed steps
$\Delta t=2\times10^{-3},2\times10^{-4},2\times10^{-5},2\times10^{-6}$ and
$2\times10^{-9}$. The minimum component changes from
$-1.54\times10^{-4}$ to $-1.49\times10^{-5}$, $-6.00\times10^{-7}$,
$+8.40\times10^{-7}$ and $+1.00\times10^{-6}$, respectively; positivity is
therefore recovered for the two smallest steps in this state. Varying the
depleted-species floor $\epsilon$ and scanning
$2,000$ logarithmically spaced steps from
$10^{-9}$ to $10^{-2}$ gives a worst component of $+2.99\times10^{-4}$ at
$\epsilon=10^{-3}$, but negative minima from $-5.80\times10^{-4}$ to
$-6.78\times10^{-4}$ for $\epsilon=10^{-4}$ to $10^{-8}$. The associated
depleted-species log jump grows from $4.61$ to $16.12$. The two-cell counterexample therefore excludes unconditional
positivity, while the step/floor scan quantifies its dependence on step size
and composition; neither is a universal timestep bound. A 21-point logarithmic
refinement of the floor scan brackets the sign transition for this frozen state
between $\epsilon=6.31\times10^{-4}$ and $7.08\times10^{-4}$, corresponding to
a depleted-species log jump between $4.95$ and $5.07$; this bracket is reported
only as a resolution result, not as a transferable positivity threshold.

For a source-only step a narrower statement remains valid. If a positive root
of $\bm Y^{n+1}=\bm Y^n+\Delta t\,\bm f(\bm Y^{n+1})$ exists and is computed
accurately, convexity and detailed balance give
$\rho[\varphi(\bm{Y}^{n+1})-\varphi(\bm{Y}^n)]\le
\rho\bm{w}^{n+1}\!\cdot(\bm{Y}^{n+1}-\bm{Y}^n)
=\Delta t\,\bm{w}^{n+1}\!\cdot\bm\omega^{n+1}\le0$.
Mass conservation follows from the source balance. This conditional argument
does not prove existence of a positive root, nor does it apply to the coupled
spatial operator without an independent invariant-domain result. The ensemble
data for implicit Euler are therefore controls on converged positive solves,
not an unconditional theorem.

Table~\ref{tab:layers} places the layer-specific conclusions side by side.

For the tested WLS reconstruction, increasing $\theta$ exchanges the two-point
entropy certificate for non-orthogonal linear exactness. At $\theta=0$ the
projected two-point rate is $-0.4083$ and the linear-field error is $0.740$;
at $\theta=1$ the error is at machine zero and the global rate is $+0.1436$.
The curves cross at $\theta^\ast=0.7398$, where the linear-field error remains
$0.193$. A common $\theta$ yields a feasible one-step interval in the two
states examined, but only as a global state-dependent control and at a measured
accuracy cost.

Patankar rescaling supplies positivity and conservation at any step size for
the PDS members examined, but it does not supply the free-energy inequality.
Proposition~\ref{prop:ascent} excludes a posteriori positive scalar relaxation
along the same ascent direction. Step rejection and recomputation can recover
positive, conservative and dissipative accepted trajectories in the tested
ensemble, at the measured cost of advancing only $0.11$--$0.29$ of the
proposed physical time (Figure~\ref{fig:members_b}); this is an algorithmic
observation, not an unconditional certificate for the original update.

The practical implication is that higher-order accuracy and the three
structural properties must be assessed separately at the face, reconstruction
and time-integration stages. The present results give an explicit
counterexample to automatic inheritance at each of the latter two stages.

This separation does not preclude a composition result for operator splitting.
If the transport and source substeps use the same potential $\Phi$, each
substep is computed in its admissible domain, and each substep separately
proves $\Delta\Phi\le0$ together with conservation and positivity, then the
composition inherits those properties by telescoping the substep changes.
The counterexamples here locate the obstruction inside a substep: the WLS
reconstruction loses the spatial face decomposition, while the Patankar update
does not supply the missing thermodynamic inequality.

\begin{table}[htbp]
\centering
\caption{Layer-specific properties. Entries state whether the corresponding
property is certified, not implied, not certified, or observed only under a
state-dependent control. The last entry is the measured fraction of proposed
physical time advanced.}
\label{tab:layers}
\scriptsize
\setlength{\tabcolsep}{0pt}
\begin{tabular}{@{}>{\raggedright\arraybackslash}p{0.12\textwidth}
>{\raggedright\arraybackslash}p{0.27\textwidth}
>{\raggedright\arraybackslash}p{0.16\textwidth}
>{\raggedright\arraybackslash}p{0.16\textwidth}
>{\raggedright\arraybackslash}p{0.29\textwidth}@{}}
\toprule
Layer & Operator & Cons. & Pos. & Entropy \\
\midrule
Face & two-point, $\bm{1}^{\mathsf T}\bm{Y}_f=1$ & certified & not implied & certified if $\lambda_{\star,f}^h\ge0$ \\
Recon. & WLS, $\theta=1$ & certified & not implied & not certified \\
Recon. & WLS, common $\theta$ & certified & one step & state-dependent global rate only \\
Time & implicit Euler, source only & certified & assumed root & certified for positive root \\
Fully discrete & two-point diffusion + backward Euler & certified & state-dependent one-step control & certified if a positive root exists and $\lambda_{\star,f}^h\ge0$ at $\bm Y^{n+1}$ \\
Time & MPRK22/43, MPDeC & certified & certified & not certified \\
Time & MPRK43, step rejection/recomputation & certified & certified & observed, at $0.11$--$0.29$ time \\
\bottomrule
\end{tabular}
\end{table}

\section{Discussion}
\label{sec:discussion}

\subsection{Why the certificates do not transfer}

Section~\ref{sec:baseline} shows why the three layer-wise certificates do not
automatically compose. The obstructions are not of the same kind, and that
difference constrains possible repairs.

Between the face closure and the reconstruction the obstruction observed here
is structural. The two-point certificate is a statement about a quadratic form
that decomposes face by face; the tested cell-centred WLS tangential correction
couples faces and that decomposition is unavailable. Global dissipation on
smooth interior states is not an equivalent certificate: the twelve
state-dependent cases are globally dissipative and still carry $166$ positive
local contributions. This finding is specific to the geometry, state family
and WLS construction tested. It does not rule out other structure-preserving
spatial designs. The scheme of \citet{jungel2023} is a two-point finite-volume
construction on an admissible mesh: it couples the face mobility and entropy
variables through a vector-valued discrete chain rule and therefore illustrates
a distinct two-point design route rather than a multipoint WLS correction. The
same distinction applies to the two-point Stefan--Maxwell construction of
\citet{cances2023}.

Between the source operator and the time integrator the obstruction is
directional. Patankar rescaling constrains where the update lands, not which
way it points, and Proposition~\ref{prop:ascent} shows that positive scalar
relaxation along the same fixed ascent direction cannot help. The tested
repair rejects the proposal and recomputes a different update, which explains
its cost in physical time. This is complementary to the Lax--Wendroff extension of
\citet{bender2025}, where a discrete entropy inequality is a hypothesis on the
scheme: the present results indicate that the hypothesis is not automatic for
detailed-balance chemistry.

\subsection{Scope and limitations}

\paragraph{Transport scope} The arithmetic-weight face criterion remains
nonnegative for the production neutral-air/carbon data and the quasi-neutral
ambipolar-air control. Its activation in this study is confined to the
synthetic diffusivity populations.

\paragraph{Discretization scope} The WLS result is an existence result for the
tested multipoint construction, and the Patankar result covers the named
members. The common-$\theta$ control is a global, state-dependent construction,
not a local limiter. A general invariant-domain condition for the coupled
diffusion--reaction system and an $h$-uniform version of \eqref{eq:trunc}
remain open.

\paragraph{Model scope} The system is isothermal, uses prescribed density, and
excludes gas--surface and thermal coupling. The WLS and Patankar components are
analysed separately rather than combined in a single reacting-flow simulation;
the two-cell calculation supplies the fully discrete separation.

\section{Conclusions}
\label{sec:conclusions}

This study examined whether conservation, positivity and thermodynamic
dissipation remain compatible as a multicomponent reacting-flow discretization
is assembled across its spatial and temporal layers. At the two-point diffusion
level, the entropy-variable jump is exactly related to the mass-fraction jump
through a log-mean Gibbs Hessian. Requiring the mass-correction term to conserve
total species mass and to be entropy neutral on the simplex tangent space then
determines a unique componentwise-positive face composition. With this choice,
the correction term drops out of the tangent-space entropy criterion, leaving
the remaining condition entirely in the Fickian part. Separate arithmetic-weight
controls remain nonnegative for the production neutral-air/carbon data and the
quasi-neutral ambipolar-air model.

This face-level certificate does not automatically survive non-orthogonal
reconstruction. The two-point entropy inequality relies on a face-wise
quadratic decomposition, whereas the weighted least-squares correction couples
a wider stencil. In the explicit counterexample, the projected two-point
operator has a global entropy rate of $-0.4083$, while the fully corrected WLS
operator gives $+0.1436$. Thus conservation and linear exactness can be retained
even after the face-wise entropy certificate is lost; entropy stability of the
face closure is not, by itself, a certificate for the reconstructed spatial
operator.

The temporal layer shows a different limitation. Patankar rescaling preserves
positivity and conservation but does not determine the sign of the free-energy
change. An explicit finite-step counterexample gives a positive and conservative
step with $\Delta\Phi>0$, and the same behaviour is found for all six
third-order members examined, on $40.0$--$56.7\%$ of the tested
detailed-balance networks. Proposition~7 identifies the underlying mechanism:
if the proposed update is already an ascent direction of the convex potential,
positive scalar rescaling of that direction cannot restore dissipation.
Positivity preservation and free-energy decay therefore require distinct
conditions.

The fully discrete comparison completes the separation. A positive face entropy
margin does not ensure that the species transport points inward on the
composition simplex, and the two-cell example shows that a backward-Euler step
can still produce a negative mass fraction. Conservation, positivity and
thermodynamic dissipation therefore constrain different stages of the
discretization and are not automatically inherited from one stage to the next.
These results provide a layer-wise framework for analysing structure
preservation in reacting-flow discretizations: conservation, positivity and
thermodynamic dissipation are not interchangeable constraints, but arise from
different discrete algebraic structures that require compatible certificates
at each layer.

\section*{Supplementary material}

The supplementary material provides the reproducibility archive description,
including the source scripts, fixed data products, figure-generation sources,
numerical checks, checksum manifests, and reproduction instructions.

\section*{Acknowledgements}

This work was supported by the China Postdoctoral Science Foundation
(Grant No. 2026M794657). Additional support was provided by the Postdoctoral
Research Funding of Hangzhou International Innovation Institute, Beihang
University (Grant No. 2025BKZ044).

\section*{Declaration of interests}
The author declares that he has no known competing financial interests or
personal relationships that could have appeared to influence the work reported
in this paper.

\appendix

\section{Proofs}
\label{app:proofs}

\paragraph{Proposition~\ref{prop:identity}}
From \eqref{eq:potential},
$\jump{w_s}=\frac{R_u}{M_s}\big(\jump{\ln Y_s}-\jump{\ln Z}\big)$.
Apply $\jump{\ln a}=\jump{a}/a^{\log}$ to both logarithms and use
$\jump{Z}=\bm{b}^{\mathsf T}\jump{\bm{Y}}$. Collecting terms gives
$\jump{\bm{w}}=\Hlog\jump{\bm{Y}}$ with $\Hlog$ as in \eqref{eq:hlog}.

\paragraph{Corollary~\ref{cor:telescope}}
Multiply the semi-discrete update by $\bm{w}_i$, sum over cells, and use
the interior-face cancellation on the closed mesh. Each interior face
contributes $|f|\jump{\bm{w}}^{\mathsf T}\bm{J}_f$, and
Proposition~\ref{prop:identity} converts this into
$-\tfrac{\rho_f|f|}{d_{n,f}}\jump{\bm{Y}}^{\mathsf T}\Hlog\bm{A}_f\jump{\bm{Y}}$.
Only the symmetric part of $\Hlog\bm{A}_f$ survives the quadratic form.

\paragraph{Proposition~\ref{prop:kernel}}
Applying \eqref{eq:hlog} to $\Ylog$ gives
$\Hlog\Ylog=R_u\big(\bm{b}-\zeta\bm{b}\big)=R_u(1-\zeta)\bm{b}$.
If $\zeta=1$ then $\Ylog$ spans the kernel, since $\Hlog$ has rank $N-1$ by the
same rank-one structure. If $\zeta\ne1$ the image of $\Ylog$ is a non-zero
multiple of $\bm{b}$, and $\Hlog$ is invertible.

\paragraph{Proposition~\ref{prop:tangentweight}}
Let $\bm D=\operatorname{diag}(1/(M_sY_s^{\log}))$, so that
$\Hlog/R_u=\bm D-\bm b\bm b^{\mathsf T}/Z^{\log}$. For $\zeta<1$,
Sherman--Morrison gives
\begin{equation}
R_u(\Hlog)^{-1}\bm 1
=\bm D^{-1}\bm 1+
\frac{\bm D^{-1}\bm b\,\bm b^{\mathsf T}\bm D^{-1}\bm 1}
{Z^{\log}-\bm b^{\mathsf T}\bm D^{-1}\bm b}
=\big(M_sY_s^{\log}\big)_s+
\frac{s_{\log}\Ylog}{Z^{\log}(1-\zeta)}.
\end{equation}
The tangent-neutral condition is
$\bm Q^{\mathsf T}\Hlog\bm Y_f=0$. Equivalently,
$\Hlog\bm Y_f$ is a scalar multiple of $\bm 1$. Hence conservation fixes
$\bm Y_f=(\Hlog)^{-1}\bm 1/(\bm 1^{\mathsf T}(\Hlog)^{-1}\bm 1)$ uniquely.
Multiplying numerator and denominator by
$Z^{\log}(1-\zeta)/s_{\log}$ gives \eqref{eq:tangentweight}; every factor is
positive. If $\zeta=1$, symmetry implies
$\operatorname{range}\Hlog=(\operatorname{span}\{\Ylog\})^\perp$.
Since $\bm1^{\mathsf T}\Ylog=s_{\log}>0$, the condition
$\Hlog\bm Y_f=c\bm1$ forces $c=0$, and
Proposition~\ref{prop:kernel} plus conservation gives
$\bm Y_f=\Ylog/s_{\log}$. This is also the limit of
\eqref{eq:tangentweight}. Substitution of the different choice
$\bm Y_f=\Ylog/s_{\log}$ and
$\Hlog\Ylog=R_u(1-\zeta)\bm b$ gives \eqref{eq:defect}.

\paragraph{Proposition~\ref{prop:shift}}
For $\bm{z}\in T$,
$\bm{z}^{\mathsf T}\Hlog(\bm{I}-\bm{Y}_f\bm{1}^{\mathsf T})\bm{z}
=\bm{z}^{\mathsf T}\Hlog\bm{z}$, so replacing $\bm{d}$ by $\bm{d}+\delta\bm{1}$
adds $\delta\,\bm{z}^{\mathsf T}\Hlog\bm{z}$ to the numerator of the Rayleigh
quotient and leaves the denominator unchanged.

\paragraph{Proposition~\ref{prop:ascent}}
Convexity gives
$\Phi(\bm{y}+\gamma\Delta\bm y)\ge\Phi(\bm{y})+\gamma\nabla\Phi(\bm{y})\cdot\Delta\bm y
\ge\Phi(\bm{y})$ for each $\gamma\ge0$ whose point lies in
$\operatorname{dom}\Phi$.

\setcounter{figure}{0}
\setcounter{table}{0}
\section{Verification populations and supporting controls}
\label{app:verification}

Three independent checks guard against implementation error rather than against
modelling error. The face-closure propositions are re-derived in exact symbolic
algebra, so that they do not rest on floating-point sampling
(Section~\ref{sec:layer1}). The decisive time-integration counterexample is
re-run in an external package that shares no code with this work
(Section~\ref{sec:layer3:thirdparty}). The time integrators are additionally
checked on three standard production--destruction benchmarks: the six
third-order members reproduce a Robertson reference to
$8.7\times10^{-11}$--$1.5\times10^{-9}$ over short times and stay positive and
conservative over long times, a Brusselator reference error decreases with step
size, and a nonlinear algal-bloom problem attains finest observed orders of
$2.932$--$2.969$. None of the three satisfies detailed balance, so none enters a free-energy statement.

Populations are frozen before use. Face states are drawn either from bounded
random intervals or from kinetic-order transport models. Reaction networks are
random, thermodynamically consistent, and constructed to satisfy detailed
balance, with molar masses
solved from the mass-balance null space so that conservation is exact by
construction; irreducibility is tested on the mass-transfer graph at the
sampled state. The $1\,440$ entries in the multi-seed population comprise
$1\,089$ distinct computed grid--state cases arranged in $480$ paired
refinement sequences; the deterministic zero-perturbation controls are repeated
only to preserve pairing.

The face-weight population nominally contains $40\,000$ random draws. One draw
has an indistinguishable left and right state and is excluded before the
criterion is evaluated, leaving the reported $39\,999$ admissible samples.
The population definition fixed the seeds, manufactured modes, perturbation
families, grids, reference tolerances and extension rule before evaluation. The
bootstrap resampling convention is stated below.

The reported confidence intervals use $4\,000$ seed-cluster bootstrap
replicates, with the perturbation seed as the resampling unit.

\begin{table}[htbp]
\centering
\caption{Numerical evidence populations. Sizes count evaluations and are not
necessarily independent statistical units. The last column records the
outcome under the stated criteria.}
\label{tab:pops}
\footnotesize
\setlength{\tabcolsep}{3pt}
\renewcommand{\arraystretch}{0.88}
\begin{tabularx}{\textwidth}{@{}l>{\raggedright\arraybackslash}p{0.29\textwidth}>{\raggedright\arraybackslash}Xc@{}}
\toprule
Layer & Population & Quantity assessed & Outcome \\
\midrule
1 & $35\,000$ faces, $N=2,\ldots,8$ & tangent-neutral identity & criterion met \\
1 & $39\,999$ random face states & five-weight ordering & criterion met \\
1 & $5\,210$ production-profile faces & arithmetic face criterion & nonnegative \\
1 & $36\,000$ quasi-neutral ambipolar faces & arithmetic face criterion & nonnegative \\
1 & refined composition front & resolution requirement & criterion met \\
2 & $55\,400$ WLS face states & per-face certificate & local certificate lost \\
2 & $600$-state search, $32$ perturbations & sign of global rate & positive rate observed \\
2 & four WLS exponents & specific to $\kappa=2$? & not specific to $\kappa=2$ \\
2 & one frozen WLS state & distance definition & criterion met \\
2 & one frozen WLS state, independent NumPy assembly & independent reassembly & independently reproduced \\
2 & curved-cosine manufactured solution & spatial accuracy criterion & \textbf{criterion not met} \\
2 & controlled lattice displacement & non-orthogonality diagnostic criterion & \textbf{criterion not met} \\
2 & $3$ grids, $4$ geometry families & order classification & criterion met \\
2 & global entropy-rate cap trajectory & positivity under cap & \textbf{criterion not met} \\
2 & $40$ seeds, $1\,440$ records & interval sufficiency & \textbf{criterion not met} \\
3 & $150$ networks, $30$ step sizes & order or rescaling & rescaling observed \\
3 & $4\,500$ implicit-Euler attempts & nonlinear solves & $92$ unsuccessful solves \\
3 & $120$ networks, six members & member coverage & free-energy increase observed \\
3 & $6$ networks, six-level refinement & step-controlled member order & \textbf{criterion not met} \\
3 & $12$ networks, step-controlled & cost of the step control & criterion met \\
3 & one network, external package & implementation independence & independently reproduced \\
cross-layer & one two-cell implicit solve & entropy $\Rightarrow$ positivity? & implication not established \\
cross-layer & step/floor scan & mechanism of positivity loss & criterion met \\
\bottomrule
\end{tabularx}
\end{table}

\begin{figure}[htbp]
  \centering
  \begin{subfigure}{0.48\textwidth}\caption{}%
    \includegraphics[width=\textwidth]{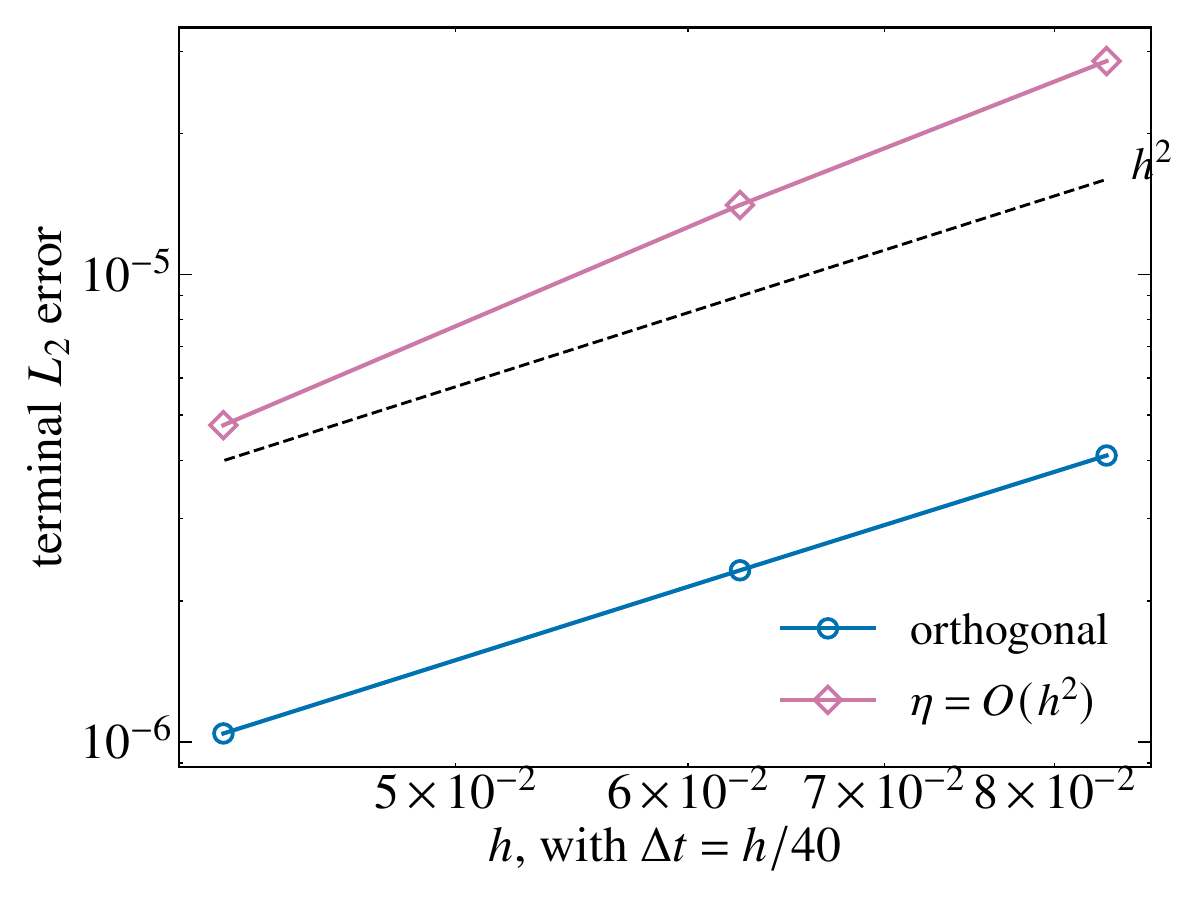}%
    % [图源] figures/scripts/fig4_pairing_and_robustness.py | 数据: data/p52_nonlinear_multicomponent_time_space_mms.json
    \label{fig:pairing_a}\end{subfigure}\hfill
  \begin{subfigure}{0.48\textwidth}\caption{}%
    \includegraphics[width=\textwidth]{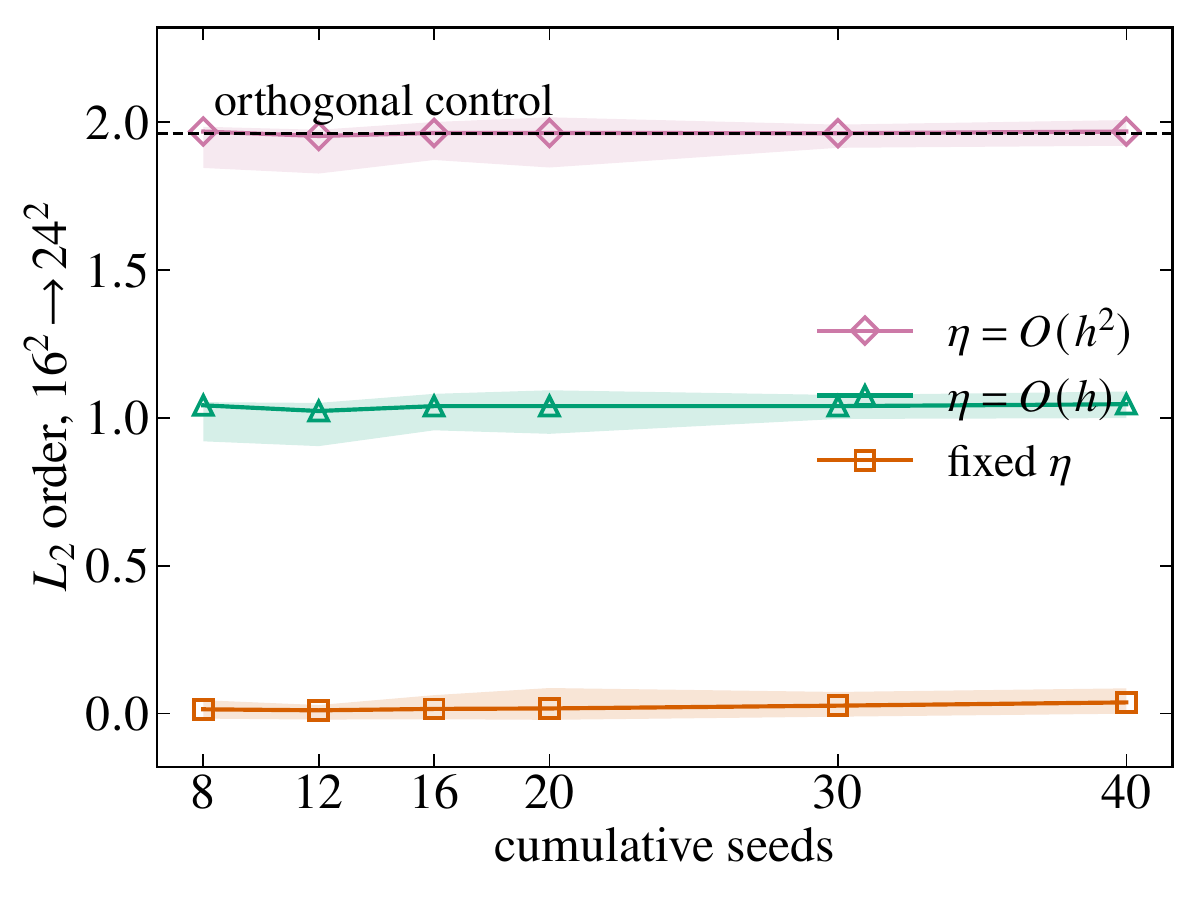}%
    % [图源] figures/scripts/fig4_pairing_and_robustness.py | 数据: data/p53_multiseed_multimode_robustness.json
    \label{fig:pairing_b}\end{subfigure}
  \caption{Paired refinement and cross-seed stability.
  (\textit{a}) Terminal $L_2$ error along the diagonal sequence
  $\Delta t=h/40$ for two geometry families.
  (\textit{b}) Cumulative seed-cluster bootstrap median and $95\%$ interval for
  the finest $L_2$ order. The orthogonal control is computed once per mode and
  grid and carries no seed variation, so it is drawn as a constant reference.
  The fixed-$\eta$ family has no convergence order to estimate, and it is that
  family whose interval fails the interval-stability target.}
  \label{fig:pairing}
\end{figure}

\section{Reproduction}
\label{app:repro}

All quantitative figures are produced by scripts under \texttt{figures/scripts/}
that read only from frozen data files. Figures~\ref{fig:certificate-map} and
~\ref{fig:geometry} are explicitly labelled analytical schematics, with
illustrative coordinates fixed in their scripts.
Each manuscript figure command (written as
\texttt{\textbackslash includegraphics}) is
annotated with its generating script and data path. The robustness
design, its seeds, modes and thresholds are hashed before use, and every data
product carries a recorded SHA-256 checksum. Every numeric value quoted in the
text is re-read from its authoritative data file by a verification script, which
distinguishes reported values from stated bounds. The complete numerical
verification record,
including criteria that were not met and the stated numerical tolerances used to assess them,
is kept with the data.

The detailed values compressed out of the main narrative remain available here:
the fixed-$\eta$ paired-refinement medians are $0.0380$ in $L_2$ and $-0.3626$
in $L_\infty$; the third-order population has largest mass error
$2.2\times10^{-12}$. Of the $4\,500$ implicit-Euler attempts, $55$ raw
unconverged iterates are retained in the numerical record as excluded
candidates; they are not counted as admissible roots. The independent Julia reproduction gives relative state
difference $1.23\times10^{-15}$, minimum component $1.13\times10^{-2}$,
MPRK22 gaps $0.0158499$ and $0.0025765$, ascent products $0.0100794$ and
$0.00253109$, and backward-Euler control gaps $-5.81\times10^{-3}$ and
$-6.09\times10^{-2}$ at $\Delta t=1$ and $100$, respectively. These are
supporting reproduction values, not additional claims.

From the project root, the complete manuscript reproduction check is run as
\begin{quote}\small\ttfamily
bash manuscript/reproduce.sh
\end{quote}
In the submitted supplementary archive, the equivalent command is run from
the archive root as \texttt{bash reproduce.sh}; the project-tree path above is
provided only for development use. The critical face, reconstruction,
time-integration and two-cell results are regenerated before the check with
\begin{quote}\small\ttfamily
bash manuscript/reproduce.sh \verb|--refresh-critical|
\end{quote}
The latter command runs the tangent-neutral identity, the frozen face-weight
population, the spatial positivity counterexample, the implicit-Euler
nonlinear-solve classification, the definition-consistent spatial reassembly,
the step/jump mechanism check and the independent NumPy-only reassembly of the
Layer-2 counterexample. The default command checks every quoted number, runs
the regression suite and compiles this PDF. The code, frozen input data and the
checksum manifest are supplied with the manuscript as supplementary material.
The default reproduction is self-contained; only the optional external
Patankar cross-check requires Julia and the external
\texttt{PositiveIntegrators.jl} package. No upstream solver source is modified
by the scripts in this archive.

\bibliographystyle{elsarticle-harv}
% The web-documentation entries contain long, unbreakable URLs; allow modest
% interword stretch in the bibliography so that the rendered references do not
% acquire underfull boxes while preserving the body-text line breaking.
\begingroup\sloppy
\bibliography{references}
\endgroup

\end{document}